\documentclass[3p]{elsarticle}

\usepackage[english]{babel}
\usepackage{amsmath}
\usepackage{amsfonts}
\usepackage{amssymb}
\usepackage{amsthm}
\usepackage{graphicx}
\usepackage{array}
\usepackage{tabularx}
\usepackage{color}
\usepackage{subfigure}
\usepackage{cancel}
\usepackage{float}
 \usepackage[normalem]{ulem}
 \newcommand{\R}{\mathbb R}

 \newcommand{\Q}{\mathbb Q}
 
\newcommand{\PP}{\mathbb{P}}

\newcommand{\dt}{\Delta t}

\newcommand{\bu}{\mathbf{u}}

\newcommand{\bxx}{{\mathbf {x}}}
\newcommand{\bn}{{\mathbf {n}}}

\newcommand{\hf}{\hat{f}}

\newcommand{\dpar}[2]{\dfrac{\partial #1}{\partial #2}}

\newcommand{\red}[1]{#1}

\definecolor{uclagold}{rgb}{1.0, 0.7, 0.0}
\definecolor{cadmiumorange}{rgb}{0.93, 0.53, 0.18}
\definecolor{darkspringgreen}{rgb}{0.18, 0.55, 0.34}
\renewcommand{\sout}[1]{}
\begin{document}
	
\title{{A high-order nonconservative approach\\ for hyperbolic equations in fluid dynamics}}

\author[uzh]{R.~Abgrall}
\ead{remi.abgrall@math.uzh.ch}

\author[uzh]{P.~Bacigaluppi}
\ead{paola.bacigaluppi@math.uzh.ch}

\author[uzh]{S.~Tokareva}
\ead{svetlana.tokareva@math.uzh.ch}

\address[uzh]{Institute of Mathematics, University of Zurich, Switzerland}

\date{\today}

\begin{abstract}
It is well known, thanks to Lax-Wendroff theorem, that the local conservation of a numerical scheme for a conservative hyperbolic system is a simple and systematic way to guarantee that, if stable, a scheme will provide a sequence of solutions that will converge to a weak solution of the continuous problem. In \cite{LeFlochHou}, it is shown that a nonconservative scheme will not provide a good solution. The question of using, nevertheless, a nonconservative formulation of the system and getting the correct solution has been a long-standing debate. In this paper, we show how to get a relevant weak solution from a pressure-based formulation of the Euler equations of fluid mechanics. This is useful when dealing with nonlinear equations of state because it is easier to compute the internal energy from the pressure than the opposite. This makes it possible to get oscillation-free solutions, contrarily to classical conservative methods. An extension to multiphase flows is also discussed.\sout{{, as well as a multidimensional extension.}}
\end{abstract}

\begin{keyword}
 Nonconservative formulation, residual distribution, conservation, fluid dynamics, {Euler equations, multiphase flow systems}
\end{keyword}

\maketitle

\section{Introduction}
According to the Lax-Wendroff theorem, it is well known that, when considering the numerical approximation of a system of hyperbolic PDEs written in conservative form, the numerical scheme must be written in conservation form, too. It is also known that, for a sequence of meshes with characteristic sizes tending to zero, if a sequence of solutions remains bounded and if its subsequence converges in some norm in $L^p$, $p\geq 1$, then the limit solution is the weak solution of the original PDE. Moreover, if the scheme satisfies a discrete entropy inequality, then the limit solution will automatically satisfy an entropy inequality. If conservation is lost, then there is no hope to get any meaningful solution, see \cite{LeFlochHou} for the analysis.

However, for engineering purposes, the conservative formulation of the behavior of a mechanical system is not necessarily the best one. Consider for example the Euler equations of fluid dynamics. The system of PDEs is
\begin{equation}
\label{eq:cons}
\dpar{}{t}\begin{pmatrix} \rho \\ \rho \mathbf{u} \\ E\end{pmatrix} +\text{ div} \begin{pmatrix}\rho \mathbf{u}\\ \rho \mathbf{u}\otimes \mathbf{u}+p\,\text{Id}\\(E+p)\mathbf{u}\end{pmatrix}=0,
\end{equation}
supplemented by initial and boundary conditions.
As usual, $\rho$ stands for the density, $\mathbf{u}$ for the velocity, and the total energy is $$E=e+\frac{1}{2}\rho \mathbf{u}^2.$$ The pressure $p$ is related to these variables via an equation of state (EOS):
\begin{equation}
\label{eos}
p=p(\rho, e)=p(\rho, E-\frac{1}{2}\rho \mathbf{u}^2).
\end{equation}
The system \eqref{eq:cons} is hyperbolic if 
$\kappa:=\dpar{p}{e}>0$
since the speed of sound $c$ is defined by
$$c^2=\kappa h, \quad h=\dfrac{e+p}{\rho}.$$
However, for engineering purposes, the relevant variables are not the conserved ones but rather the primitive ones, namely density, velocity and internal energy or pressure. When the solution is smooth, system \eqref{eq:cons} can be equivalently written as:
\begin{equation}
\label{eq:nc}
\dpar{}{t}\begin{pmatrix}\rho\\ \rho \mathbf{u} \\ e\end{pmatrix} +\begin{pmatrix}\text{ div } \rho \mathbf{u}\\
\text{div }(\rho \mathbf{u}\otimes \mathbf{u}+p\,\text{Id})\\
\mathbf{u}\cdot \nabla e+(e+p) \text{ div }\mathbf{u} \end{pmatrix}
=0
\end{equation}
or
\begin{equation}
\label{eq:nc:bis}
\dpar{}{t}\begin{pmatrix}\rho\\ \rho \mathbf{u} \\ p\end{pmatrix} +\begin{pmatrix}\text{ div } \rho \mathbf{u}\\
\text{div }(\rho \mathbf{u}\otimes \mathbf{u}+p\,\text{Id})\\
\mathbf{u}\cdot \nabla p+\rho c^2 \text{ div }\mathbf{u} \end{pmatrix}
=0.
\end{equation}
These equations are valid for smooth flows and cannot be considered for discontinuities. Nevertheless, there have been several attempts to solve the Euler equations {with either } formulation \eqref{eq:nc} or \eqref{eq:nc:bis}, including solutions with shocks. One example of  such method{s} is Karni's hybrid scheme \cite{karni1,karni2}, where formulation \eqref{eq:nc:bis} is used along slip lines only thanks to a switch in the scheme. In any case, this method violates strict conservation.

This has been a long ongoing debate on how to, nevertheless, use formulations \eqref{eq:nc} or \eqref{eq:nc:bis} for the numerical approximation of Euler equations valid for all kinds of flows, e.~g. with complex equations of state. In the case of nonlinear equations of state, i.~e. when the pressure explicitly depends on the density, the pressure obtained by the numerical scheme cannot be uniform across contact discontinuities. The reason for that behaviour is that, on one hand, the density is evaluated from the mass conservation, and, on the other hand, one evaluates the pressure via energy and density. If, in addition, we want the pressure to be constant across contact discontinuities, this puts a constraint that is in general not compatible with the updated densities, momentum and total energy, see \cite{abgrallhandbook} for a short discussion.

{In} our knowledge, in the Eulerian framework, there exists only one approach to this problem which is described in \cite{saurel}. In this paper, we propose a simpler and more general framework for dealing with nonconservative formulations. We solve equations \eqref{eq:nc} or \eqref{eq:nc:bis} in a way which is compatible with local conservation and the continuity of pressure and velocity across contact discontinuities. To achieve this, we rely on a finite volume formulation that uses residuals instead of fluxes. In a flux formulation, the unknowns are approximations of the average values of the conserved variables, and they are balanced by the sum of normal fluxes across the boundary of the control volume. This assumes that the control volume has a polygonal shape. In general, these control volumes are interpreted as cells of a dual mesh which is made of simplices. In the residual formulation, one starts by a mesh whose elements are simplices, and interprets the unknowns as approximations of the point values of the conserved variables. These unknowns, for any given degree of freedom, are then updated by a sum of the local residuals over all the elements that share this degree of freedom. Given any element $K$, the local conservation is recovered by requiring that the sum of the local residuals for that element is the normal flux over the boundary of $K$ of some consistent approximation of the flux. It is easy to show that any flux formulation leads to a residual form, and the opposite is also true. However, the fluxes that are computed depend not only on the solution on both sides of the face of the control volume, see \cite{icm} for details.

The format of this paper is the following. In Section~\ref{From conservative to non conservative formulations}, we recall how one can get a residual distribution formulation for the system \eqref{eq:nc} that is equivalent to a flux formulation of \eqref{eq:cons}. This enables us to get a relation on the increment of the energy that can be generalised for the residual formulation. In Section~\ref{Conservative approximation}, we show how to use this principle, first on {the} energy-based formulation of the Euler equations \eqref{eq:nc} and then on the pressure-based formulation \eqref{eq:nc:bis} for several kinds of equations of state. In Section~\ref{multiphase}, this is further generalised to multiphase flows where the phases may have very complex equations of state. Finally, we give some concluding remarks in Section~\ref{Conclusions}.

\section{From {a} conservative to {a} nonconservative formulation}\label{From conservative to non conservative formulations}

\subsection{A residual formulation of a finite volume scheme}

The main advantage of the residual formulation can be understood from the one dimensional setting. Consider the problem:
$$\dpar{U}{t}+\dpar{f(U)}{x}=0.$$
With standard notations, a generic finite volume {scheme} writes:
$$U_j^{n+1}=U_j^n-\lambda\big (\hf_{j+1/2}-\hf_{j-1/2})$$
where $\lambda=\Delta t/\Delta x$ and $\hf_{j+1/2}$ is the flux between the states $U_{j}^n$ and $U_{j+1}^n$. High order accuracy in space amounts to tune the arguments of the flux, {while high order in time can be reached via a strong stability preserving (SSP) scheme.} To fix the conservation problem one must fix the scheme to recover a flux form, i.e to work directly with the fluxes. This is not easy from an algebraic point of view. 

It is known, see for example \cite{abgrall2001} that any finite volume scheme can be rewritten in terms of a distribution of the residual. Consider for example, and for simplicity, the one dimensional case, its generalisation to any kind of control volume is straightforward, see again \cite{abgrall2001}.

The residual formulation writes (in its simplest form) as
$$U_j^{n+1}=U_j^{n}-\lambda(\Phi_j^{j+1/2}+\Phi_j^{j-1/2}).$$ The conservation is recovered if for any element $[x_j, x_{j+1}]$ one gets
\begin{equation}\label{eq:cons:rd}
\Phi_j^{j+1/2}+\Phi_{j+1}^{j+1/2}=f(U_{j+1})-f(U_j)
\end{equation}
for any order of accuracy.
{\sout{For example one} One} can go from a flux formulation to a residual formulation by defining{, for example}:
$$\Phi_j^{j+1/2}=\hf_{j+1/2}-f(U_j) {\quad} \text{and} {\quad} \Phi_{j+1}^{j+1/2}=f(U_{j+1})-\hf_{j+1/2}.$$
The  local conservation is a consequence of relation \eqref{eq:cons:rd}. If we start from a nonconservative formulation in residual form, one can check the conservation if one can provide linear transformations of {this} residual to obtain a form satisfying \eqref{eq:cons:rd}.

\subsection{Unsteady residual distribution formulation for the conservative case}\label{unsteady:cons}
Consider a multidimensional hyperbolic system in the form
\begin{equation*}
\dpar{U}{t}+\text{div }\mathbf{f}(U)=0.
\end{equation*}

Recall the residual distribution approach from \cite{ricchiuto}. We start by a Runge-Kutta (RK) formulation: knowing $U^n$, we define:
\begin{equation*}
\begin{split}
U^{(0)}&=U^n,\\
\frac{U^{(1)}-U^{{(0)}}}{\Delta t}&+\text{div } \mathbf{f}(U^{(0)})=0,\\
\frac{U^{(2)}-U^{{(0)}}}{\Delta t}&+\frac{1}{2}\bigg ( \text{div } \mathbf{f}(U^{(0)})+\text{div } \mathbf{f}(U^{(1)})\bigg )=0,\\
U^{n+1}&=U^{(2)}.
\end{split}
\end{equation*}
Next for $l=0,1$ we rewrite each sub-step as:
\begin{equation}
\label{rk:rd}
\frac{U^{(l+1)}-U^{{(0)}}}{\Delta t}+ \text{ DIV }\mathcal{F}(U^{(l)},U^{{(0)}})=0
\end{equation}
where
\begin{equation}
\label{rkrd:2}
\text{ DIV }\mathcal{F}(U^{(l)},U^{{(0)}})=\frac{1}{2}\bigg (  \text{div } \mathbf{f}(U^{(0)})+\text{div } \mathbf{f}(U^{(l)})\bigg )
\end{equation}
This will provide an approximation that is second order in time.
Without loss of accuracy, we can also write it as 
\begin{equation}
\label{rkrd:3}
\text{ DIV }\mathcal{F}(U^{(l)},U^{{(0)}})=\text{div } \mathbf{f}\bigg (\frac{  U^{(0)}+U^{(l)}}{2}\bigg )
\end{equation}
and the adapted modifications for the RK scheme. 

We assume that the computational domain $\Omega$ is covered by non-overlapping simplices $\{K_j\}_{j\in \mathcal{J}}$, $\Omega=\cup_{j\in \mathcal{J}}\{K_j\}$. The elements $K_j$ are segments in 1D,  triangles/quadrilateral in 2D and tetrahedrons/hexahedrons in 3D. In order to simplify the notations, we denote by $\PP^1(K)$ the set of polynomials of degree $1$ on triangles/tetrahedrons or by $\Q^1{(K)}$ the set {of} polynomials of degree $1$ on quadrilaterals/hexahedrons. Both guarantee a second order accurate approximation of any $C^1$ function. 
We introduce the approximation space
$$V_h=\{ u\in L^2(\Omega) : \text{ for any }K_j, j\in \mathcal{J}, u_{|K_j}\in \PP^1(K_j)\}.$$
We denote by {$\sigma$ the degrees of freedom}, i.e. the vertices of the element $K_j$. 
Then the residual distribution approximation of \eqref{rk:rd} reads: for any degree of freedom $\sigma$, define first $U^{{(0)}}=U^n$ and for $l=0,1$, do
\begin{equation}
\label{rk:rd:approx}
\begin{split}
|C_i| \dfrac{ U_\sigma^{(l+1)}-U_\sigma^{{({0})}}}{\Delta t}&+\sum_{K\ni \sigma} \Phi_\sigma^K(U^{(l)}_h,U^{{(0)}}_h)=0\\
{\text{with}\quad}  \Phi_\sigma^K(U^{(l)}_h,U^{{(0)}}_h)&=\beta_\sigma^{K}(U^{(l)}_h,U^{{(0)}}_h) \bigg ( \int_{K} \dfrac{U_h^{(l)}-U_h^{{(0)}}}{\dt} dx + \int_{\partial K}\mathcal{F}(U^{(l)}_h,U^{{(0)}}_h)\cdot \bn {\,d\Gamma} \bigg )
\end{split}
\end{equation}
In \eqref{rk:rd:approx}, the parameter $ \beta_\sigma^{K}(U^{(l)}_h,U^{{(0)}}_h)$ is a matrix ({for Euler system of equations} $3\times 3$ in 1D, $4\times 4$ in 2D, etc.), and it is constructed as a formal extension of the same construction for scalar problems that is guaranteed to have an $L^\infty$ stability bound, see \cite{ricchiuto} for details. 

To compute $\beta_\sigma$, we  can proceed as follows.
\begin{itemize}
\item In \eqref{rk:rd:approx}, instead of the $\beta$-formulation, we first consider the Rusanov residual:
\begin{equation}
\label{Rusavov:xt}
\Phi_\sigma^{K,Rus}(U^{(l)}_h,U^{{(0)}}_h)=\frac{1}{N_K}\bigg ( \int_{K} \dfrac{U_h^{(l)}-U_h^{{(0)}} }{\dt} dx + \int_{\partial K}\mathcal{F}(U^{(l)}_h,U^{{(0)}}_h)\cdot \bn{\,d\Gamma}\bigg )+\alpha_K \big ( \frac{U_h^{(l)}+U_h^{{(0)}}}{2}-\widehat{U}_K\big )
\end{equation}
with
$$\widehat{U}_K=\sum_{\sigma\in K} \frac{U_\sigma^{(l)}+U_\sigma^{{(0)}}}{2}.$$
In \eqref{Rusavov:xt}, $N_K$ corresponds to the number of degrees of freedom in $K$.

\item Then we consider the eigen-decomposition of $A_\bn=A(\widehat{U}_K)n_x+B(\widehat{U}_K)n_y$ (in 2D, while in 1D it would simply be $A(\widehat{U}_K)$. The matrices $A$ and $B$ are the Jacobians of the $x$- and $y$- component of the flux $\mathbf{f}$ with respect to the state $\widehat{U}_K$. Here, the vector $\bn$ is a unit vector in the direction of the velocity field when it is nonzero, or any arbitrary direction otherwise. Of course, this direction is not relevant in one dimension. The right eigenvectors of $A_\bn$ are denoted by
$\{\mathbf{r}_\xi\}$. We denote by $\{\mathbf{\ell}_\xi\}$ the left eigenvectors of $A_\bn$, so that any state $X$ can be written as:
$$X=\sum_\xi \mathbf{\ell}_\xi(X)\mathbf{r}_\xi.$$
%
%

\item Then we decompose, for any degree of freedom $\sigma$:
$${\Phi_\sigma^{K,Rus}(U^{(l)}_h,U^{{(0)}}_h) =\sum_{\xi} \mathbf{\ell}_{\xi}\big [ \Phi_\sigma^{K,Rus}(U^{(l)}_h,U^{{(0)}}_h)\big ] \mathbf{r}_{\xi}}$$
We note that for any $\xi$,
$$\mathbf{\ell}_\xi\big [ \Phi^{K}\big ]=\sum_{\sigma \in K} \mathbf{\ell}_\xi\big [ \Phi_\sigma^{K,Rus}(U^{(l)}_h,U^{{(0)}}_h)\big ],$$
where $\Phi^{K} = \sum_{\sigma \in K}\Phi_\sigma^{K,Rus}$ is the total residual.

\item The next step is to define $\Phi_\sigma^{K,\star}(U^{(l)}_h,U^{{(0)}}_h)$ as:
\begin{subequations}\label{grosbeta}
\begin{equation}
\label{star}
\Phi_\sigma^{K,\star}(U^{(l)}_h,U^{{(0)}}_h)=\sum_\xi \beta_\sigma^{K,\star} \mathbf{\ell}_\xi\big ( \Phi^{K}\big )\mathbf{r}_\xi
\end{equation}
where:
\begin{equation}
\label{betastar}
\beta_\sigma^{K,\star}=(1-\Theta_\xi)
\dfrac{ 
\max \Big( \frac{\mathbf{\ell}_\xi\big ( \Phi_\sigma^{K,Rus}(U^{(l)}_h,U^{{(0)}}_h)\big )}{\mathbf{\ell}_\xi\big [ \Phi^{K}\big ]},0\Big)
}
{\sum\limits_{\sigma'\in K}\max \Big( \frac{\mathbf{\ell}_\xi\big ( \Phi_{\sigma'}^{K,Rus}(U^{(l)}_h,U^{{(0)}}_h)\big )}{\mathbf{\ell}_\xi\big [ \Phi^{K}\big ]},0\Big)
}
+\Theta_\xi \mathbf{\ell}_\xi\big ( \Phi_\sigma^{K,Rus}(U^{(l)}_h,U^{{(0)}}_h)\big )
\end{equation}
with
\begin{equation}
\label{theta}
\Theta_\xi=\dfrac{\big | \mathbf{\ell}_\xi\big (\Phi_\sigma^{K,Rus}(U^{(l)}_h,U^{{(0)}}_h)\big )\big |}{\sum\limits_{\sigma'\in K}\big | \mathbf{\ell}_\xi\big ( \Phi_{\sigma'}^{K,Rus}(U^{(l)}_h,U^{{(0)}}_h)\big )\big |}.
\end{equation}
\end{subequations}
Note that $\Theta_\xi\in [0,1]$.
This guarantees that the scheme is second order in time and space and (formally) non-oscillatory, see \cite{ricchiuto,Abgrall2006} for more details. 
\end{itemize}

\subsection{Discrete conservation for a nonconservative formulation}\label{cons:nc}

As a first exercise, let us show that any explicit scheme for \eqref{eq:cons} can be rewritten as an explicit scheme for \eqref{eq:nc} or \eqref{eq:nc:bis}.
It is well known that
$$dE=de +\mathbf{u}\cdot d\mathbf{m} -\frac{1}{2}\;||\mathbf{u}||^2 d\rho,$$ 
where $\mathbf{m}=\rho\mathbf{u}$ is the momentum, and thus the question is to see how one can use this relation for discrete problems.

To simplify the notations, we start from a simple Euler time stepping method, {but it is important to note, that \sout{while}} more general algorithms can be treated similarly.
Setting $U=(\rho, \mathbf{m}, E)^T$, {where $\mathbf{m}=\rho\mathbf{u}$,} we start from the finite volume scheme
\begin{equation}
\label{eq:fv}
|C_i| \big ( U_i^{n+1}-U_i^n)+\sum_{j\in \mathcal{V}_i} \hat{\mathbf{f}}_{ij}=0
\end{equation}
where $C_i$ is the area/volume of the control volume, $\mathcal{V}_i$ is the set of neighbouring cells to cell $i$, and $\hat{\mathbf{f}}_{ij}$ is the numerical flux between the cells $C_i$ and $C_j$. In the following, since $\hat{\mathbf{f}}_{ij}$ doesn't play any role, we write:
$$\delta \hat{\mathbf{f}}:=\sum_{j\in \mathcal{V}_i} \hat{\mathbf{f}}_{ij}.$$

This finite volume scheme in the component-wise form reads:
\begin{subequations}
\label{euler:cons}
\begin{align}
\label{euler:cons:rho}
|C_i|(\rho_i^{n+1}-\rho_i^n)+\delta\hat{f}_\rho &=0,\\
\label{euler:cons:m}
|C_i|(\mathbf{m}_i^{n+1}-\mathbf{m}_i^n)+\delta\hat{f}_{\mathbf{m}} &=0,\\
\label{euler:cons:E}
|C_i|(E_i^{n+1}-E_i^n)+\delta\hat{f}_E &=0.
\end{align}
\end{subequations}

Next, we introduce the $\Delta$ operator
\begin{subequations}\label{Delta}
\begin{equation}\label{Delta:01}
\Delta g=g^{n+1}-g^n
\end{equation}
and make the observation that
\begin{equation}\label{Delta:02}
\Delta (gh)= g^{n+1}\Delta h+h^n\Delta g,
\end{equation}
which we can rewrite as 
\begin{equation}
\label{Delta:1}
\Delta (gh)=\overline{g}\Delta h+\underline{h}\Delta g
\end{equation}
with
\begin{equation}
\label{Delta:2}
\overline{g}=g^{n+1}, \qquad \underline{g}=g^n,
\end{equation}
so that
\begin{equation}
\label{Delta:3}
\overline{gh}=\overline{g}\; \overline{h}\qquad \underline{gh}=\underline{g}\; \underline{h}
\end{equation}
\end{subequations}
From \eqref{Delta} we see that
\begin{equation*}
\Delta E=\Delta e+\frac{1}{2}\Delta (\rho \mathbf{u}^2)=\Delta e+\frac{1}{2}\bigg( \overline{\rho \mathbf{u}}\cdot \Delta \mathbf{u}+\underline{\mathbf{u}}\cdot \Delta \mathbf{m}\bigg)
\end{equation*}
and since $\overline{\rho \mathbf{u}}=\overline{\rho}\; \overline{ \mathbf{u}}$ and
$\Delta \mathbf{m}=\overline{\rho}\Delta \mathbf{u} +\underline{\mathbf{u}}\Delta \rho$, we can write
\begin{equation}
\Delta E=\Delta e+\frac{\overline{\mathbf{u}}+\underline{\mathbf{u}} }{2}\cdot \Delta \mathbf{m}-\frac{1}{2}\overline{\mathbf{u}}\cdot \underline{\mathbf{u}}\; \Delta \rho.
\end{equation}
From this last equation we see that 
\begin{subequations}
\label{euler:nc:e}
\begin{equation}
\label{euler:cons:e}
|C_i|(e_j^{n+1}-e_j^n)+ \delta \hat{f}_e=0
\end{equation}
with
\begin{equation}
\label{euler:cons:e2} \delta \hat{f}_e=\delta \hat{f}_E-\frac{\overline{\mathbf{u}}+\underline{\mathbf{u}} }{2}\cdot \delta \hat{f}_{\mathbf{m}}+\frac{1}{2}\overline{\bu}\cdot \underline{\bu}\; \delta \hat{f}_\rho=0\, ,
\end{equation}
{which coupled to \eqref{euler:cons:rho} and \eqref{euler:cons:m}} provides a consistent discretisation of \eqref{eq:nc:bis} that is equivalent to the discretisation \eqref{euler:cons} of \eqref{eq:cons}. 
\end{subequations}
The main problem is that { \eqref{euler:cons:rho}-\eqref{euler:cons:m}-\eqref{euler:nc:e} \sout{\eqref{euler:nc:e}-\eqref{euler:cons:m}-\eqref{euler:cons:rho}}} is identical to \eqref{euler:cons} with the same properties and same drawbacks, so no progress has yet been made.

Let us again consider the scheme \eqref{eq:fv}, and in particular its flux term.  It is known, see for example \cite{abgrall2001}, that any finite volume scheme can be rewritten in terms of distribution of a residual. Consider for simplicity, the one dimensional case, its generalisation to any kind of control volume is straightforward, see again \cite{abgrall2001}.
Relation \eqref{eq:fv}, using standard notations, writes:
$$
|C_j| (U_j^{n+1}-U_j^n)+\Delta t \big (\hat{f}_{j+1/2}-\hat{f}_{j-1/2}\big )=0,$$ i.e.
\begin{equation}\label{rd}
|C_j| (U_j^{n+1}-U_j^n)+\Delta t \big (\Phi_j^{j+1/2}+\Phi_j^{j-1/2}\big )=0\, ,
\end{equation}
where we have set 
\begin{equation*}
\begin{split}
\Phi_j^{j+1/2}&=\hat{f}_{j+1/2}-f(U_j)\\
\Phi_j^{j-1/2}&=f(U_j)-\hat{f}_{j-1/2}\, .
\end{split}
\end{equation*}
Hence, on the element $[x_j,x_{j+1}]$ we can define two residuals
\begin{equation}
\label{rd:2}
\begin{split}
\Phi_j^{j+1/2}&=\hat{f}_{j+1/2}-f(U_j)\\
\Phi_{j+1}^{j+1/2}&=f(U_{j+1})-\hat{f}_{j+1/2}
\end{split}
\end{equation}
and we notice the conservation relation:
\begin{equation}
\label{rd:3}
f(U_{j+1})-f(U_j)=\Phi_{j}^{j+1/2}+\Phi_{j+1}^{j+1/2}.
\end{equation}

Coming back to the system \eqref{euler:cons}, written as a residual distribution scheme in the form  \eqref{rd}, the residuals $\Phi_\sigma^{j+1/2}$, $\sigma=j$, $j+1$, have 3 components, namely {$\Phi_{\sigma,\rho}^{j+1/2}$, $\Phi_{\sigma,\mathbf{m}}^{j+1/2}$ and $\Phi_{\sigma,E}^{j+1/2}$} for the mass, momentum and total energy, respectively. From \eqref{euler:cons:e} and \eqref{euler:cons:e2} we see that we can define a residual for the internal energy as follows:
\begin{equation*}
{\Phi_{\sigma,e}^{j+1/2}=\Phi_{\sigma,E}^{j+1/2}-\frac{\overline{\mathbf{u}}_\sigma+
\underline{\mathbf{u}}_\sigma }{2}\cdot \Phi_{\sigma,\mathbf{m}}^{j+1/2}+\frac{1}{2}\overline{\bu}\cdot \underline{\bu}\;\cdot \Phi_{\sigma,\rho}^{j+1/2}.}
\end{equation*}
We obviously have
\begin{equation}
\label{cons:e}
{\sum_{\sigma=j}^{j+1}\Phi_{\sigma,e}^{j+1/2}=\sum_{\sigma=j}^{j+1}\Phi_{\sigma,E}^{j+1/2}- \sum_{\sigma=j}^{j+1}\frac{\overline{\mathbf{u}}_\sigma+\underline{\mathbf{u}}_\sigma }{2}\cdot \Phi_{\sigma,\mathbf{m}}^{j+1/2}+\frac{1}{2}\sum_{\sigma=j}^{j+1}\overline{\bu}\cdot \underline{\bu}\;\cdot \Phi_{\sigma,\rho}^{j+1/2}.}
\end{equation}
Defining the total residuals as
$${\Phi_\xi^{j+1/2}=\sum_{\sigma=j}^{j+1}\Phi_{\sigma,\xi}^{j+1/2}}$$
for $\xi=e$ and $E$, we see that \eqref{cons:e} rewrites as
\begin{equation}
\label{cons:e total}
{\Phi_e^{j+1/2}=\Phi_E^{j+1/2}-\sum_{\sigma=j}^{j+1}\frac{\overline{\mathbf{u}}_\sigma
+\underline{\mathbf{u}}_\sigma }{2}\cdot \Phi_{\sigma,\mathbf{m}}^{j+1/2}+\frac{1}{2}\sum_{\sigma=j}^{j+1}\overline{\bu}\cdot \underline{\bu}\;\cdot \Phi_{\sigma,\rho}^{j+1/2}.}
\end{equation}
Relation \eqref{cons:e total} represents our target relation.

\section{Conservative approximation of the Euler equations in primitive variables}
\label{Conservative approximation}

If we are able to define the residuals ${\Phi_{j,e}^{j+1/2}}$ and ${\Phi_{j+1,e}^{j+1/2}}$ that  satisfy \eqref{cons:e}-\eqref{cons:e total}, then using the same technique as in \cite{Abgrall:Roe}, we can show that if the conditions of the Lax-Wendroff theorem  for the sequence of approximation hold, then the limit solution is a weak solution of \eqref{eq:cons}. A sketch of the proof is recalled in \ref{appendix A}. In this section, we show how to achieve this requirement on system \eqref{eq:nc}  using the residual formulation described in Section \ref{unsteady:cons}. To describe the principle of the method, we start with system \eqref{eq:nc} and  we deal then later with system \eqref{eq:nc:bis} in the case of non linear equations of state. 
The discussion will be general, however, the numerical experiments will be one dimensional in order to compare the results with exact solutions. 

\subsection{Residual distribution formulation for unsteady problems in non-conservative form}\label{unsteady:nc}

To describe our method we first consider system \eqref{eq:nc} and start by considering the first order case, and then extend it to the second order in time and space.

\subsubsection{First order scheme} 
We start by {a} first order case in order to illustrate the method. The temporal scheme is a simple one step Euler scheme, and the space residual will also be first order accurate. We set $U=(\rho, \mathbf{m}, e)^T$, and knowing $U^n$ at every degree of freedom, $U^{n+1}$ is obtained by:
$$|C_\sigma|(U_\sigma^{n+1}-U_\sigma^n)+\dt \sum_{K\ni\sigma}\Phi_\sigma^K(U^n)=0.$$
This scheme has the same format as the one of Section \ref{From conservative to non conservative formulations} that has been used to derive \eqref{cons:e total}. 
To be specific, we calculate the Rusanov residual for the first two components of $U$ by applying formula \eqref{Rusavov:xt}, which for the first order gives
\begin{equation}
\begin{split}
&\Phi_{\sigma,\rho}^{K,Rus}=\frac{1}{N_K} \bigg ( \int_{\partial K}\mathcal{F}_\rho(U^n)\cdot \bn{ \, d\Gamma }\bigg )+\alpha_K \big ( \rho_{\sigma}^n-\widehat{\rho}_{K}\big ),\\[0.5em]
&\Phi_{\sigma,\mathbf{m}}^{K,Rus}=\frac{1}{N_K} \bigg ( \int_{\partial K}\mathcal{F}_\mathbf{m}(U^n)\cdot \bn{ \, d\Gamma }\bigg )+\alpha_K \big ( \mathbf{m}_\sigma^n-\widehat{\mathbf{m}}_K\big ).
\end{split}
\end{equation}

For the internal energy, we define the Rusanov residual as:
\begin{equation}
\label{eq:Phi_Rus}
\Phi_{\sigma,e}^{K,Rus}=\frac{1}{N_K} \int_K \big (\bu\cdot \nabla e+\rho h \text{ div }\bu\big )^{\star\star}\,d\bxx + \alpha_K\big ( {e_\sigma^n}-\widehat{e}_K\big ) \, ,
\end{equation}
The choice of the approximation $(\bu\cdot \nabla e+\rho h \text{ div }\bu\big )^{\star\star}$ in \eqref{eq:Phi_Rus} is dictated only by the accuracy constraint, namely that for a smooth solution we must have
$$\int_K \big (\bu\cdot \nabla e + \rho h \text{ div }\bu\big )^{\star\star}\,d\bxx-\int_K \big (\bu \cdot \nabla e + \rho h\text{ div } \bu \big ) d\bxx=O(h^{d+k}),$$
where $k$ is the degree of the interpolation and $d$ is the number of spacial dimensions,
see \cite{abgrall:eleve} for more details.

Experimentally, one can see that under a CFL like constraint, the numerical solution converges, but the limit is not a weak solution of the problem, as expected. The reason is that the conservation relation \eqref{cons:e} is not satisfied.  Hence, from the density and the momentum equation, we can compute the velocity at time $t_{n+1}$, and then we modify the residual on the internal energy by {setting}:
\begin{equation}
\label{energy:perturb}
{\widetilde{\Phi_{\sigma,e}^{K,Rus}}} := \Phi_{\sigma,e}^{K,Rus}+r_\sigma^K.
\end{equation}
The perturbations $r_\sigma^K$ are chosen such that the relation \eqref{cons:e} is satisfied,
\begin{equation*}
\sum\limits_{\sigma\in K}\Phi_{\sigma,e}^{K,Rus}+\sum\limits_{\sigma\in K}r_\sigma^K={\Phi_E^{K}}- \sum\limits_{\sigma\in K}\frac{\bu_\sigma^{n+1}+\bu_\sigma^n}{2}\cdot 
 \Phi_{\sigma, \mathbf{m}}^{K,Rus}+\frac{1}{2}\sum\limits_{\sigma\in K}\bu_\sigma^n\cdot \bu_\sigma^{n+1}\;{\cdot \,}\Phi_{\sigma,\rho}^{K,Rus}.
\end{equation*}
There is no reason to favour one degree of freedom over another, therefore we set:
\begin{equation}
\label{energy:perturb:2}
\begin{split}
r_\sigma^K=r^K&=\frac{1}{N_K} \bigg ({\Phi_E^{K}}- \sum\limits_{\sigma\in K}\frac{\bu_\sigma^{n+1}+\bu_\sigma^n}{2}\cdot \Phi_{\sigma, \mathbf{m}}^{K,Rus}\\
&\qquad \qquad\qquad+\frac{1}{2}\sum\limits_{\sigma\in K}\bu_\sigma^n\cdot \bu_\sigma^{n+1}\;{\cdot \,}\Phi_{\sigma,\rho}^{K,Rus}-\sum\limits_{\sigma\in K}\Phi_{\sigma,e}^{K,Rus}\bigg )\, .
\end{split}
\end{equation}

%

\subsubsection{Second order scheme}\label{second order case}
We shall use scheme \eqref{rk:rd} to get second order of accuracy. One of its properties is that again each stage of the scheme looks like a first order Runge-Kutta scheme but with a modified residual since one needs to take into account both the time and space increments in relation
\eqref{cons:e total}.  More specifically, we set $U=(\rho, \mathbf{m}, e)^T$ and proceed as follows.
\begin{itemize}
\item For the first step: we first compute the velocity $\bu^{(1)}$ and  relation \eqref{cons:e total} is written with
$$\overline{\bu}=\bu^{(1)}, \qquad \underline{\bu}=\bu^{{(0)}}=\bu^n,$$
\item For the second step, knowing the states $U^{{(0)}}=U^n$ and $U^{(1)}$, we first compute $\bu^{(2)}=\bu^{n+1}$ and then
$$\overline{\bu}=\bu^{(2)}, \qquad \underline{\bu}=\bu^{(1)}.$$
These modifications are mandatory to get proper conservation { (see  \eqref{Delta:01}-\eqref{Delta:02})}.
\end{itemize}
This means that ${\Phi_E^K}$ is now:
$${\Phi_E^K}=\int_K \dfrac{E^{(l)}-E^{{(0)}}}{\Delta t} {d \mathbf{x}}+\int_{\partial K} \mathcal{F}_E(U^{(l)},U^{{(0)}})\cdot \bn \, {d \Gamma \,},
$$
where the total energy flux $\mathcal{F}_E$ is evaluated according to \eqref{rkrd:2} or \eqref{rkrd:3}. The approximation $(\bu\cdot \nabla e+\rho h \text{ div }\bu\big )^{\star\star}$ in \eqref{eq:Phi_Rus} satisfies second order accuracy e.g. by means of arithmetic averages of the variables $U^{(0)}$ and $U^{(l)}$ with $l=0,1$. The correction is:
\begin{equation}
\label{energy:perturb:2.2}
\begin{split}
r_\sigma^K=r^K&=\frac{1}{N_K} \bigg ( {\Phi_E^K}- \sum\limits_{\sigma\in K}\frac{\bu_\sigma^{\red{(l+1)}}+\bu_\sigma^{\red{(l)}}}{2}\cdot \Phi_{\sigma, \mathbf{m}}^{K,Rus}\\
&\qquad \qquad\qquad+\frac{1}{2}\sum\limits_{\sigma\in K}\bu_\sigma^{\red{(l+1)}}\cdot \bu_\sigma^{\red{(l)}}\;{\cdot \,}\Phi_{\sigma,\rho}^{K,Rus}-\sum\limits_{\sigma\in K}\Phi_{\sigma,e}^{K,Rus}\bigg )\, .
\end{split}
\end{equation}

{Remark: Once again it is important to emphasize, that it is possible to use the velocity quantities at $l+1$ to evaluate \eqref{energy:perturb:2.2}, since one first solves the mass and momentum equations for $l+1$ and only then focuses on the nonconservative internal energy (or pressure) equation.}

\subsection{Ensuring conservation in the case of non-linear equations of state}\label{3.2}
In the introduction, we have mentioned the problem of evaluating the pressure across contact discontinuities for non linear equations of state. Consider an equation of state (EOS) to be
$$p=p(\rho, e)$$
and assume it to be such that $$
\kappa= \dpar{p}{e}>0.$$ Note that for a perfect gas $\gamma=\kappa+1$, see \cite{menikoff} for details. Let us recall an example from \cite{saurel}.
We consider Cochran and Chan EOS given by
\begin{subequations}
\label{cochran}
\begin{equation}
\label{cochran:1}
p=\Gamma \rho \big (\varepsilon-\varepsilon_0(\rho)\big ) +p_0(\rho)\, ,
\end{equation}
with
\begin{equation}
\label{cochran:2}
\begin{split}
\varepsilon_0(\rho)&= \dfrac{A_1}{\rho_0(E_1-1)}\bigg (\dfrac{\rho}{\rho_0}\bigg )^{E_1-1}- \dfrac{A_2}{\rho_0(E_2-1)}\bigg (\dfrac{\rho}{\rho_0}\bigg )^{E_2-1}\, ,\\
p_0(\rho)&=A_1\bigg (\dfrac{\rho}{\rho_0}\bigg )^{E_1}-A_2 \bigg (\dfrac{\rho}{\rho_0}\bigg )^{E_2}\, .
\end{split}
\end{equation}
\end{subequations}
We set up a shock tube problem with the EOS parameters listed in Table \ref{tab:cochran} and consider piecewise-constant initial data given by
\begin{itemize}
\item $u=1000$ [m/s], $p=20\cdot10^9$ [Pa]
\item $\rho_L=1134$, $\rho_R=500$ [kg/m$^3$]
\end{itemize}
\begin{table}[h]
\begin{center}
\begin{tabular}{|cccccc|}\hline
$\rho_0$ & $A_1$ & $E_1$& $A_2$& $E_2$ & $\Gamma$\\
(kg/m$^3$) & (GPa) &         & (GPa) &            &                   \\
\hline
$1134$ & $0.819181$ & $4.52969$ & $1.50835$ & $1.42144$ & $1.19$\\
\hline
\end{tabular}
\caption{\label{tab:cochran} Parameters for the Cochran \& Chan EOS \eqref{cochran}.}
\end{center}
\end{table}

\begin{figure}[H]
\subfigure[velocity]{\includegraphics[width=0.45\textwidth]{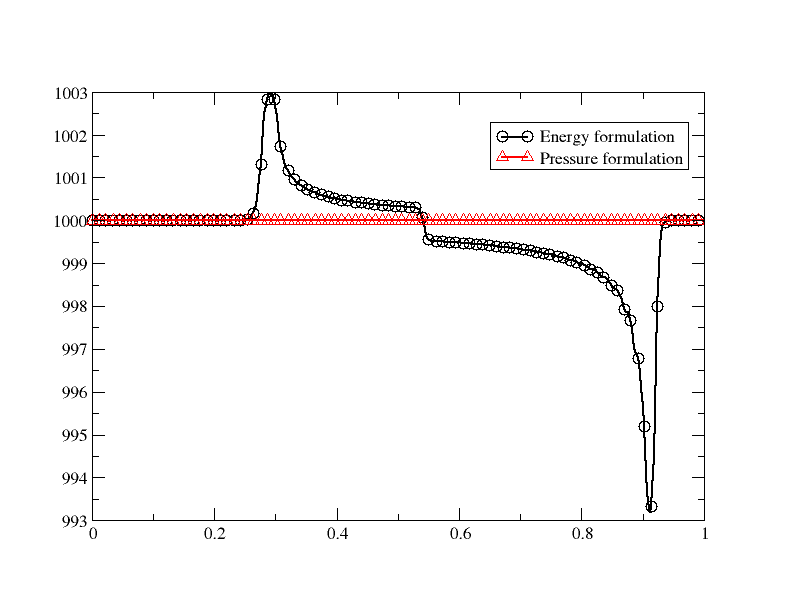} } \subfigure[velocity, taken from \cite{saurel}]{\includegraphics[width=0.45\textwidth]{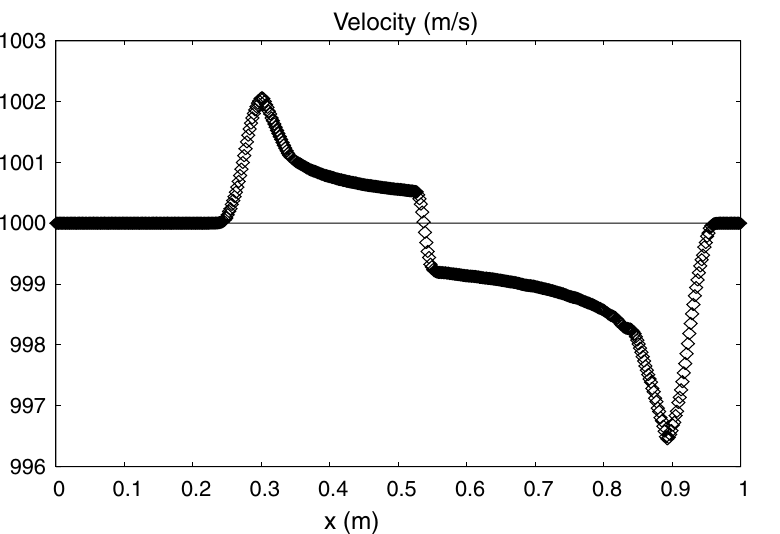} } \\
\subfigure[pressure]{\includegraphics[width=0.45\textwidth]{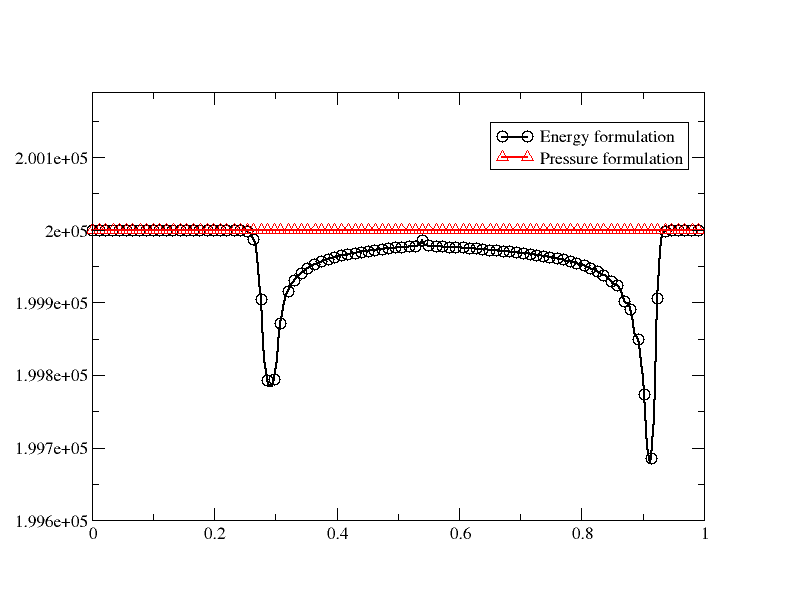} }\subfigure[pressure, taken from \cite{saurel}]{\includegraphics[width=0.45\textwidth]{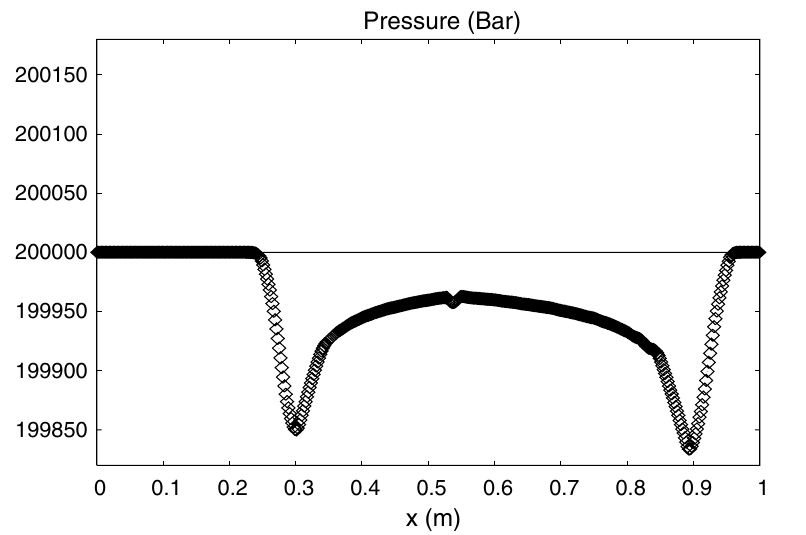} }\\
{\subfigure[density]{\includegraphics[width=0.45\textwidth]{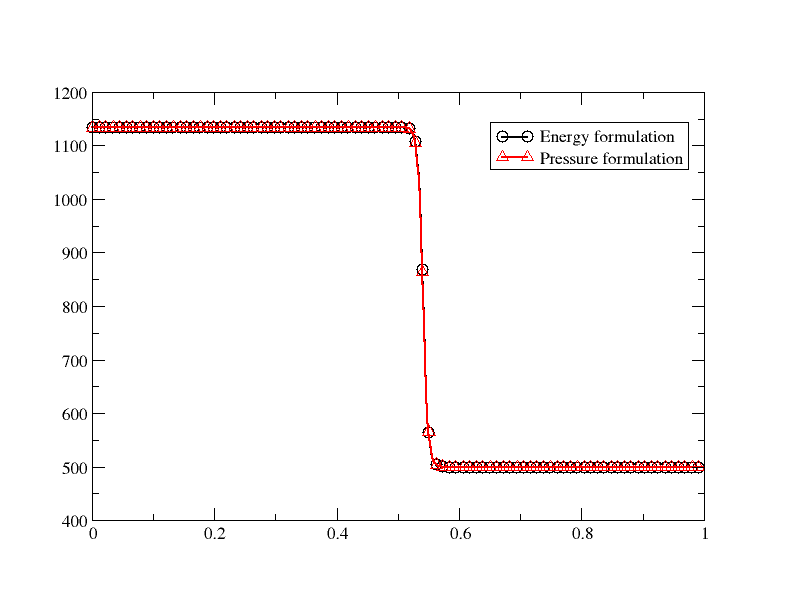} }}{\subfigure[density, taken from \cite{saurel}]{\includegraphics[width=0.45\textwidth]{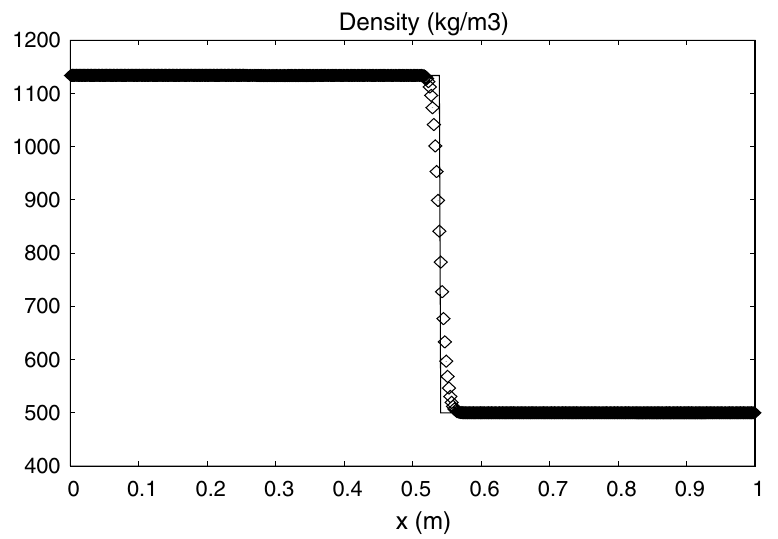} }}
\caption{\label{fig:cochran} 
Results for a contact discontinuity for Cochran and Chan EOS. The energy and pressure formulations are considered on the left. The results are compared with those obtained for an HLLC scheme of second order (taken from \cite{saurel}) on the right column.}
\end{figure}

Using the conservative scheme described in Section \ref{Conservative approximation} we get the results displayed in Figure~\ref{fig:cochran}. This results are compared to those of \cite{saurel}, 
which have been obtained with a HLLC scheme with a MUSCL extrapolation on the physical variables.  
In both cases, we see that the pressure and the velocity do not stay uniform as they should be.
It is well known that even with the Godunov scheme, the evaluation of the pressure across the contact discontinuity can be problematic, even for single fluids, see \cite{abgrallhandbook,toro-anomalies-2002} for the analysis.
The reason is that the equation of state that relates the internal energy to the density and the pressure can be highly nonlinear. The internal energy is obtained from the total energy and the kinetic energy, and thus the pressure via:
$$e=\langle E\rangle -\frac{1}{2} \dfrac{ \langle \rho \mathbf{u}\rangle ^2}{\langle \rho \rangle}
$$
where we have put \emph{conserved variables}  emphasised by using the notation $\langle \cdot \rangle$.
The problem is that {across \sout{in}} a contact, the pressure is uniform, and in the case of  a highly nonlinear EOS, there is no reason that the relation 
$$\langle E\rangle -\frac{1}{2} \dfrac{ \langle \rho \mathbf{u}\rangle ^2}{\langle \rho \rangle}= e\big (\langle \rho\rangle, p\big )$$
will guarantee a uniform pressure when $\langle \rho \rangle$ changes across a contact discontinuity. Up to our knowledge, the only Eulerian method that provides correct values is described in \cite{saurel}. It is however quite complicated and tuned for Cartesian meshes.
  
A way to solve this issue is to start from \eqref{eq:nc:bis} and use the same ideas as before. 
At the continous level, we have 
\begin{equation}\label{de}de=\dpar{e}{\rho} d\rho+\dpar{e}{p}dp,\end{equation}
therefore
$$ dE= \dpar{e}{\rho} d\rho +\dpar{e}{p} dp +\mathbf{u}\cdot d(\mathbf{m})-\frac{ \mathbf{u}^2 }{2} d\rho,$$ 
and this is the relation to mimic in the numerical scheme.
Hence, we wish to satisfy
\begin{equation}\label{constraint}
\begin{split}
  \sum_{\sigma \in K} \bigg ( \widetilde{\dpar{e}{\rho}}_\sigma \Phi_{\sigma,\rho}^K&+\widetilde{\dpar{e}{p}}_\sigma\big ( \Phi_{\sigma,p}^K+r_\sigma^p\big )\bigg ) +
  \sum_{\sigma\in K}\overline{\mathbf{u}}_\sigma\cdot \Phi_{\sigma,\mathbf{m}}^K -\sum_{\sigma\in K} \frac{\overline{\mathbf{u}^2}_\sigma}{2}\; {\cdot \,}\Phi_{\sigma,\rho}^K=
  {\Phi_E^K}\\&=\int_K \big (E^{(l)}-E^{{{(0)}}}\big ){d \mathbf{x} }+ \Delta t \int_{\partial K} \mathcal{F}_E(U^{(l)},U^{{{(0)}}})\cdot \bn\, {d\Gamma \,},
  \end{split}
\end{equation}
where $\widetilde{\dpar{e}{\rho}}_\sigma$, $\widetilde{\dpar{e}{p}}_\sigma$ are approximations of $\dpar{e}{\rho}$ and $\dpar{e}{p}$ at the degree of freedom $\sigma$ which need to be determined and $r_\sigma^p$ corresponds to the corrections on the pressure residuals. As before, we will assume that $r_\sigma^p=r^p$ is independent of $\sigma$. Note that we only perturb the pressure residual and not the density residual. The reason is that we wish to keep an explicit scheme: first we update the density, then the velocity and finally the pressure. If we can omit this constraint, more freedom can be obtained.

{Relation \eqref{constraint} can be rewritten as:
  \begin{equation}
   \label{modif:E2}
   \begin{split}
  \sum_{\sigma \in K} \widetilde{\dpar{e}{p}}_\sigma r_\sigma^p&={\Phi_E^{K}}-\bigg [
  \sum_{\sigma\in K}\overline{\mathbf{u}}_\sigma\cdot \Phi_{\sigma,\mathbf{m}}^K -\sum_{\sigma\in K} \frac{\overline{\mathbf{u}^2}_\sigma}{2}\; {\cdot \,}\Phi_{\sigma,\rho}^K+\sum_{\sigma \in K} \bigg( {\widetilde{\dpar{e}{\rho}}}_\sigma \Phi_{\sigma,\rho}^K + \widetilde{\dpar{e}{p}}_\sigma \Phi_{\sigma,p}^K \bigg) \bigg]\\
  &:=\Delta_e\, .
   \end{split}
  \end{equation}
}
  
The  question is how to define the terms $\widetilde{\dpar{e}{\rho}}_\sigma$  and $\widetilde{\dpar{e}{p}}_\sigma$. Once this is done we can get $r^p$. One constraint is that, if initially the pressure and the velocity  are uniform, we keep this property at the next time step. We will start to see how relation \eqref{constraint} behaves in case of a 1D flow with uniform velocities and pressures. We start by the first order scheme, and then go to the second order.
  
The goal is to find 'good' approximations of partial derivatives of the internal energy \eqref{de} so that the correction for the pressure vanishes for uniform pressures, without violating conservation. For this, we consider the behaviour of the internal energy increment. 
For any $\sigma$, we can write, for any $\lambda \in \R$
\begin{equation*}
e(p_\sigma^{(l+1)}, \rho_\sigma^{(l+1)}) -e(p_\sigma^{{{(l)}}}, \rho_\sigma^{{{(l)}}}) = \widetilde{\dpar{e}{p}}_{\sigma,\lambda}\big ( p_\sigma^{(l+1)}-p_\sigma^{(l)}\big ) + \widetilde{\dpar{e}{\rho}}_{\sigma,\lambda} \big (\rho_\sigma^{(l+1)}-\rho_\sigma^{{(l)}}\big)\, ,
\end{equation*}
with
\begin{equation}
\begin{split}
\widetilde{\dpar{e}{p}}_{\sigma,\lambda }=&\lambda \dfrac{e(p_\sigma^{(l+1)}, \rho_\sigma^{(l+1)}) -e(p_\sigma^{{(l)}}, \rho_\sigma^{(l+1)})}{p_\sigma^{(l+1)}-p_\sigma^{{(l)}}}+(1-\lambda)
\dfrac{ e(p_\sigma^{(l+1)}, \rho_\sigma^{{(l)}})-e(p_\sigma^{(l)},\rho_\sigma^{(l)})}{p_\sigma^{(l+1)}-p_\sigma^{(l)}}\, ,\\
\widetilde{\dpar{e}{\rho}}_{\sigma,\lambda }=&\lambda \dfrac{ e(p_\sigma^{{(l)}}, \rho_\sigma^{(l+1)})-e(p_\sigma^{(l)},\rho_\sigma^{(l)})}{\rho_\sigma^{(l+1)}-\rho_\sigma^{(l)}}+(1-\lambda)\dfrac{e(p_\sigma^{(l+1)}, \rho_\sigma^{(l+1)}) -e(p_\sigma^{(l+1)}, \rho_\sigma^{{(l)}})}{\rho_\sigma^{(l+1)}-\rho_\sigma^{{(l)}}}\, .
\end{split}
\label{tilde_def}
\end{equation}
When $p$ and $u$ are uniform, then the right hand side of \eqref{constraint} reduces to
$$	{\Phi_E^{K}}=\int_K (e^{(l)}-e^{{(0)}} ){d \mathbf{x}}+\Delta t \int_{\partial K}\dfrac{e ^{(\textcolor{black}{l})}    \mathbf{u}^{(\textcolor{black}{l})}+e^{(0)}\mathbf{u}^{(0)}}{2}\cdot \bn\,{d \Gamma\,} .$$

If for any $\lambda$ we have $\widetilde{\dpar{e}{p}}_{\sigma,\lambda }=0$, then $\Delta_e=0$.
Otherwise, one can find    $\lambda\in \R$  such that if 
$$ \dfrac{e(p_\sigma^{(l+1)}, \rho_\sigma^{(l+1)}) -e(p_\sigma^{{(l)}}, \rho_\sigma^{(l+1)})}{p_\sigma^{(l+1)}-p_\sigma^{{(l)}}} \neq 
\dfrac{ e(p_\sigma^{(l+1)}, \rho_\sigma^{{(l)}})-e(p_\sigma^{(l)},\rho_\sigma^{(l)})}{p_\sigma^{(l+1)}-p_\sigma^{(l)}}\, ,$$
then 
\begin{subequations}
\begin{equation}
\label{condition:1}
\sum_{\sigma \in K} \widetilde{\dpar{e}{p}}_{\sigma,\lambda} r_{\sigma}^p=\Delta_e
\end{equation} 
and since $\widetilde{\dpar{e}{p}}_{\sigma,\lambda } \textcolor{black}{\neq} 0$, we can find $r_{\sigma}^p=r^p\in \R$ such that
\begin{equation}
\label{condition:p}
r^p \bigg (\sum_\sigma \widetilde{\dpar{e}{p}}_{\sigma,\lambda } \bigg ) =\Phi^{E,K}-\bigg [
  \sum_{\sigma\in K}\overline{\mathbf{u}}_\sigma\cdot \Phi_{\sigma,\mathbf{m}}^K -\sum_{\sigma\in K} \frac{\overline{\mathbf{u}^2}_\sigma}{2}\; \Phi_{\sigma,\rho}^K+\sum_{\sigma \in K} \bigg( {\widetilde{\dpar{e}{\rho}}}_{\sigma,\lambda } \Phi_{\sigma,\rho}^K + \widetilde{\dpar{e}{p}}_{\sigma,\lambda } \Phi_{\sigma,p}^K \bigg) \bigg]\, .
\end{equation}

Clearly, we have to set $r^p=0$, if $p$ is uniform so that uniformity is preserved at the next time step. 
\end{subequations}

\vspace{0.5cm}
\textcolor{black}{It is important to note, that \eqref{tilde_def} is explicit in pressure in the case that the chosen equation of state guarantees that 
\begin{equation}
\label{explic_eos}
\frac{\partial p}{\partial \rho}=\psi(\rho), \quad
\frac{\partial p}{\partial e}=\eta(\rho).
\end{equation} 
This is for example always achieved when using the Mie-Gr\"uneisen EOS
$$p=\Gamma \rho(\varepsilon-\varepsilon_0(\rho))+p_0(\rho),$$ as done in the current work with the Cochran Chan EOS.
In case the reader might be interested in an EOS not fulfilling \eqref{explic_eos}, it is possible to simplify \eqref{tilde_def} and consider
\begin{equation}
\label{simple} 
\widetilde{\dpar{e}{p}}_{\sigma,\lambda }=\dpar{e}{p}(p^{(l)}_{\sigma},\rho^{(l)}_{\sigma}), \quad \widetilde{\dpar{e}{\rho} }_{\sigma,\lambda }=\dpar{e}{\rho}(p^{(l)}_{\sigma},\rho^{(l)}_{\sigma}).
\end{equation}
}


To detect contact discontinuities and preserve the uniformity of pressure and velocity, in the one-dimensional case we proceed as follows:
\begin{itemize}
\item if 
\begin{equation}
\label{simple:2}\max\bigg (\dfrac{\big |\max\limits_{\sigma\in K}\overline{u}_\sigma -\min\limits_{\sigma\in K}\overline{u}_\sigma \big |}{\big |\max\limits_{\sigma\in K}\overline{u}_\sigma \big |+\big |\min\limits_{\sigma\in K}\overline{u}_\sigma \big |+\varepsilon_1}, \dfrac{\big |\max\limits_{\sigma\in K}p^{(l)}_\sigma -\min\limits_{\sigma\in K}p^{(l)}_\sigma \big |}{\big |\max\limits_{\sigma\in K}p^{(l)}_\sigma \big |+\big |\min\limits_{\sigma\in K}p^{(l)}_\sigma \big |+\varepsilon_1}\bigg )\leq \varepsilon\, ,
\end{equation}
then we set $r^p=0$, as this would mean that we are going across a contact wave.
\item In any other case, $r^p$ is evaluated from \eqref{condition:p}.
\end{itemize}

\subsection{Numerical results}\label{Numerica illustrations}
System \eqref{eq:cons} has been tested on a set of very demanding benchmark problems with two different typologies of equations of state.

\subsubsection{Perfect gas EOS}
The first test consists in a shock tube on a domain $x=[0,1]m$ with a diaphragm located at $x_d=0.5m$ at the initial conditions left and right of $x_d$ 
\begin{itemize}
\item $\rho_l=100 \;[\frac{kg}{m^3}]$, $u_l=0 \;[\frac{m}{s}]$, $p_l=10^9 \;[Pa]$,
\item $\rho_r=1 \;[\frac{kg}{m^3}]$, $u_r=0 \;[\frac{m}{s}]$, $p_r=10^5\;[Pa]$.
\end{itemize}
The closing equations of state for system \eqref{eq:cons} are for a perfect gas and read $p(\rho,e)=(\gamma-1)\rho e$, with $\gamma=1.4$.\\

\noindent The results are shown at $t=45\; \mu s$ and have been obtained with a grid of $N=5000$ nodes and a $\text{CFL}=0.5$.
The choice of a high number of cells for this test case is intended to show how the proposed numerical approximation converges to the exact solution for a very strong "Sod"-like problem.
\begin{figure}[H]
\centerline{\subfigure[Pressure]{\includegraphics[width=0.45\textwidth]{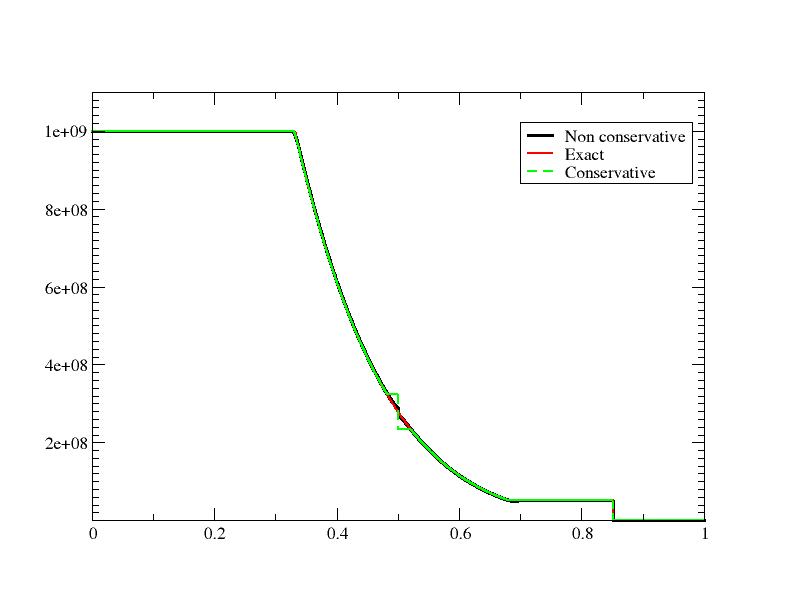}}\subfigure[Pressure, zoom]{\includegraphics[width=0.45\textwidth]{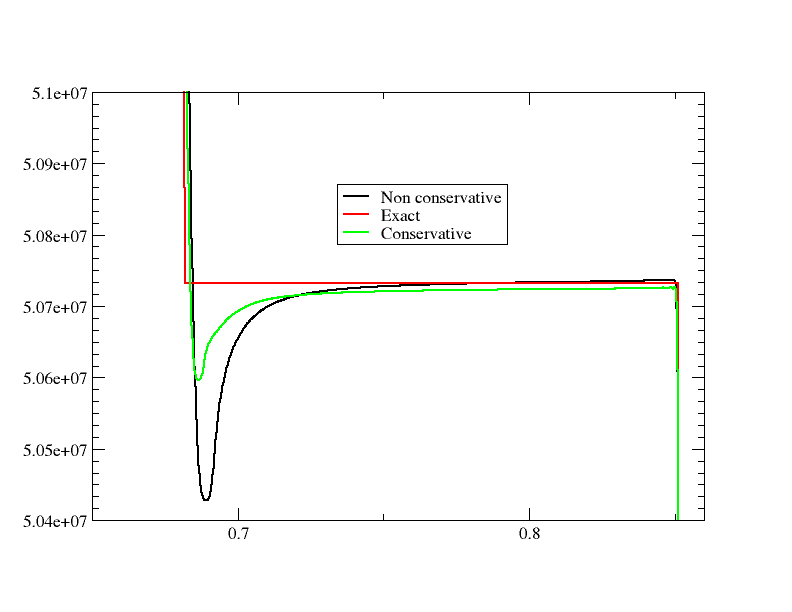}}}
\caption{\label{fig:pressureI} Strong shock tube problem with perfect gas EOS for the Euler equations. Comparison with the exact solution and both the conservative and nonconservative versions of the scheme.
 Display of the (a) pressure distribution along the shock tube and (b) a zoom on the pressure between $x=[0.65,0.86]m$ of (a).}
\end{figure}

\begin{figure}[H]
\subfigure[Density]{\includegraphics[width=0.45\textwidth]{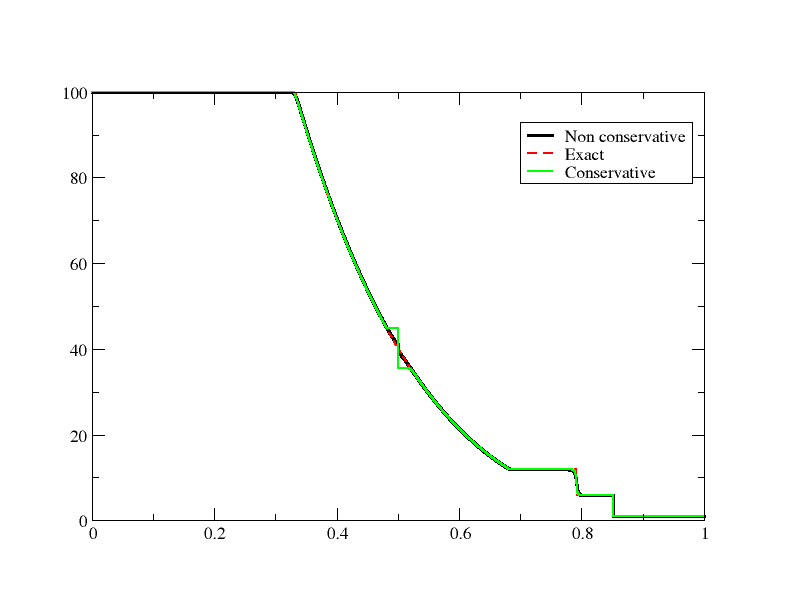}}\subfigure[Density, zoom]{\includegraphics[width=0.45\textwidth]{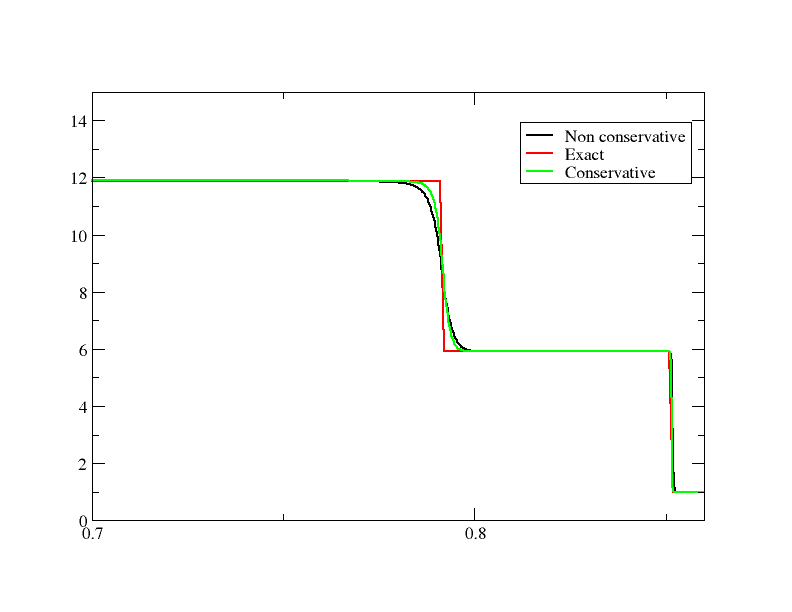}}

\subfigure[Velocity]{\includegraphics[width=0.45\textwidth]{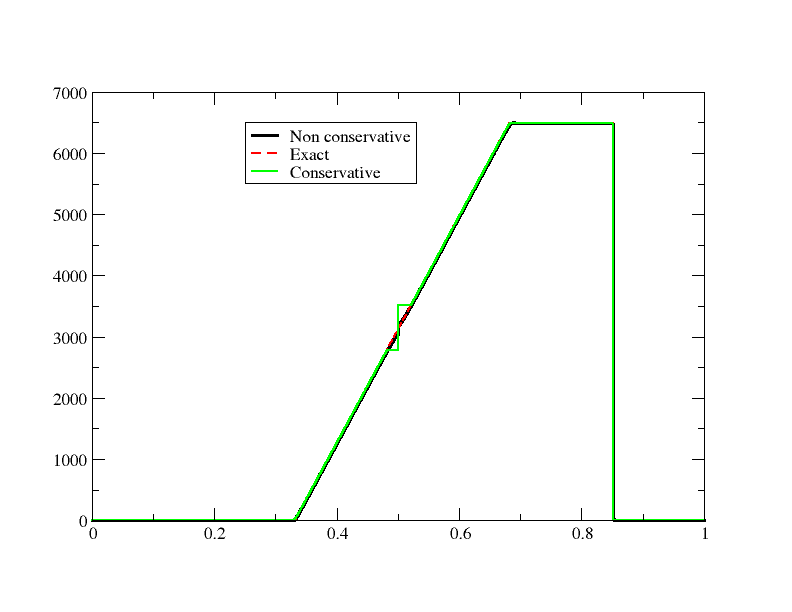}}\subfigure[{Velocity, zoom}]{\includegraphics[width=0.45\textwidth]{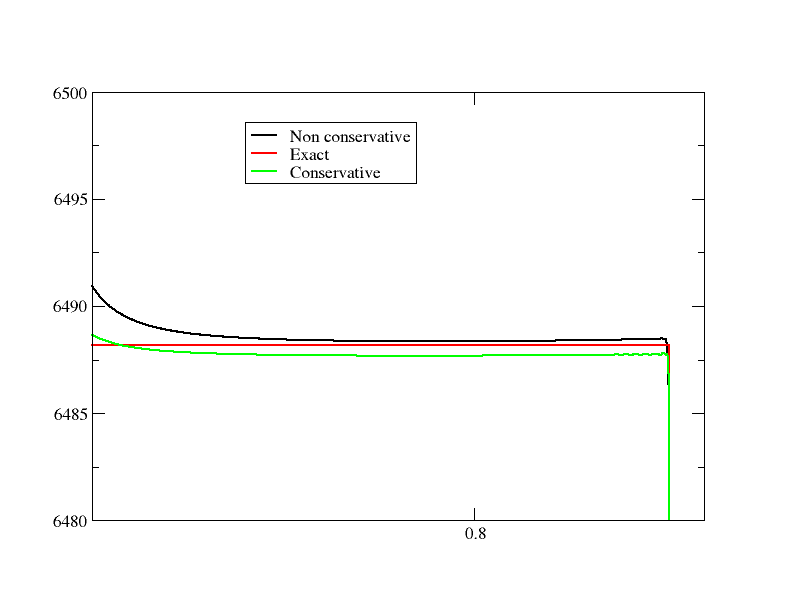}}
\caption{Strong shock tube problem with perfect gas EOS for the Euler equations. Comparison with the exact solution and both the conservative and nonconservative versions of the scheme.
Display of (a) the density distribution along the shock tube, (b) a zoom on the density  between $x=[0.65,0.86]m$, (c) the velocity distribution and (d) a zoom on the velocity between $x=[0.7,0.86]m$.\label{fig:veldensI}}
\end{figure}

Figures \ref{fig:pressureI} and \ref{fig:veldensI} show the comparison of the results given by the exact solution and the conservative and nonconservative approximations presented in this work. 
The behaviour of these solutions show a good overlap. Both the conservative and nonconservative results are characterized by a glitch at the sonic point, which is less pronounced in the nonconservative case. The glitch itself can be easily corrected by adding some entropy correction,  see e.~g.  \cite{sermeus} for a possible fix, but what matters in these results is that, this difference is also due to the fact, that the nonconservative approximation results in a slightly more diffused solution. This can be particularly seen in the zooms of the pressure (Fig. \ref{fig:pressureI}(c)), the density and of the velocity
(Fig. \ref{fig:veldensI}(b) and (d)).

\subsubsection{Nonlinear EOS}
In order to check the quality of the approximation of the nonconservative scheme both for the pressure and energy formulations, a second testcase on a Riemann problem with a strong discontinuity has been evaluated with the choice of the Cochran-Chan EOS, as described in Section \ref{3.2}.
The values for the EOS are those of Table \ref{tab:cochran} and we set\footnote{Note that the results are not sensitive to the choice of $\varepsilon$.} $\varepsilon_1=\varepsilon=10^{-6}$ in \eqref{simple:2}. The considered domain $x=[0,1]m$ is split at $x_d=0.5m$ in a left and right state, where initially the values are set to
\begin{itemize}
\item  $\rho_l=1134\;[\frac{kg}{m^3}]$, $u_l=0\;[\frac{m}{s}]$, $p_l=2\cdot 10^{10}[Pa]$
\item $\rho_r=120\;[\frac{kg}{m^3}]$, $u_r=0\;[\frac{m}{s}]$, $p_r=2\cdot 10^5[Pa]\, .$
\end{itemize}

\begin{figure}[H]
\begin{center}
\subfigure[Velocity]{\includegraphics[width=0.45\textwidth]{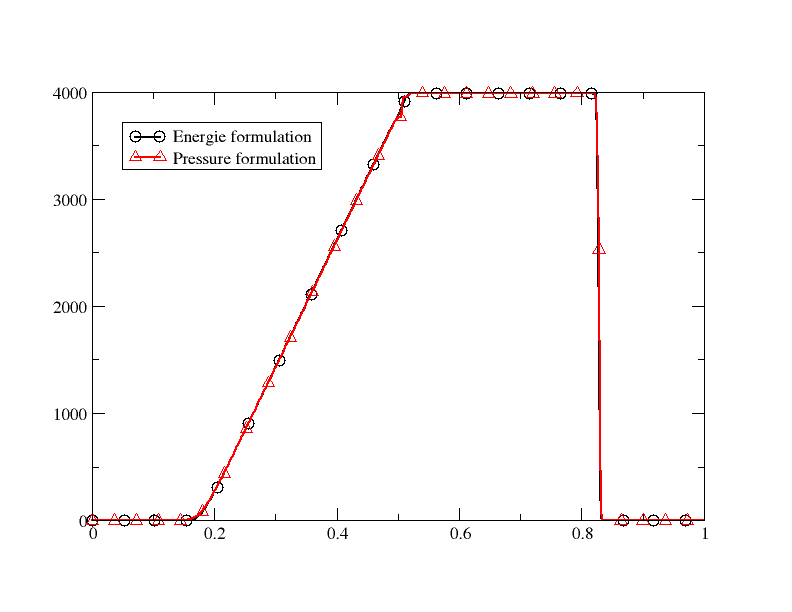}} \subfigure[Velocity, zoom]{\includegraphics[width=0.45\textwidth]{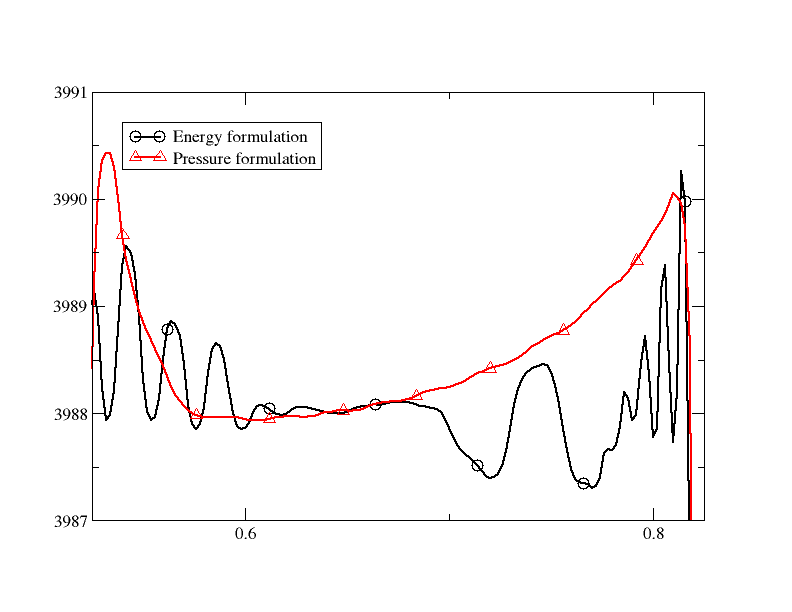}}\\
\subfigure[Density]{\includegraphics[width=0.45\textwidth]{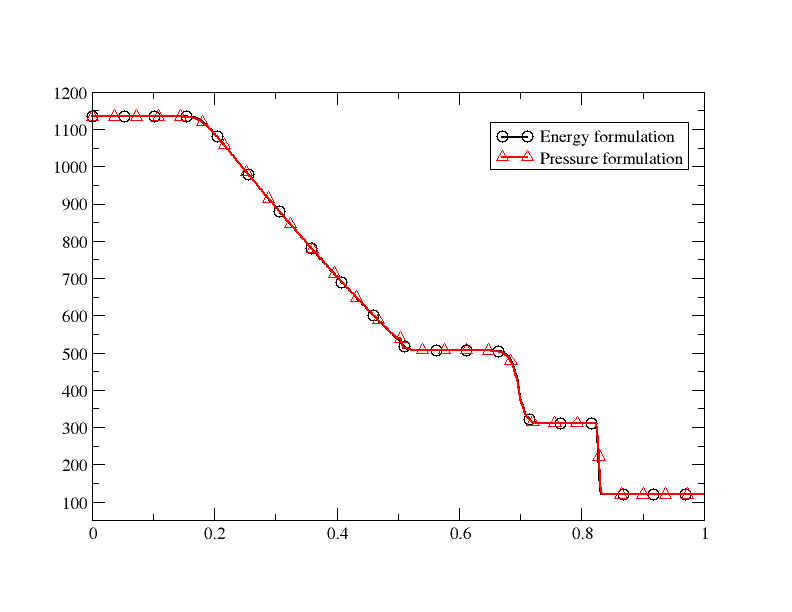}} \subfigure[Density, zoom]{\includegraphics[width=0.45\textwidth]{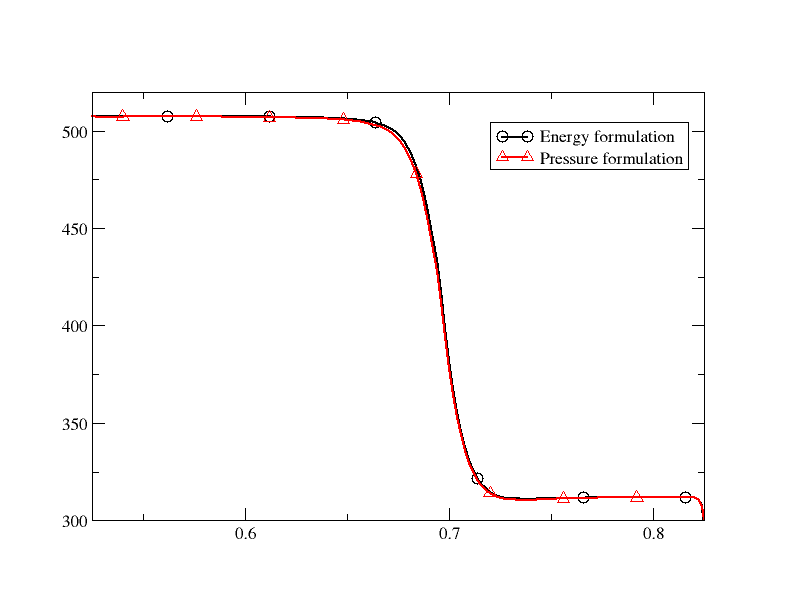}}\\
\subfigure[Pressure]{\includegraphics[width=0.45\textwidth]{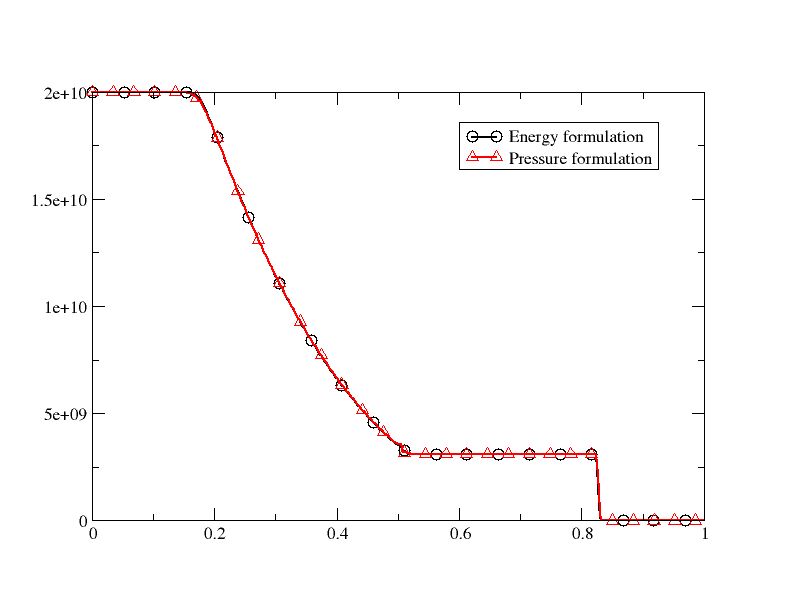}} \subfigure[Pressure, zoom]{\includegraphics[width=0.45\textwidth]{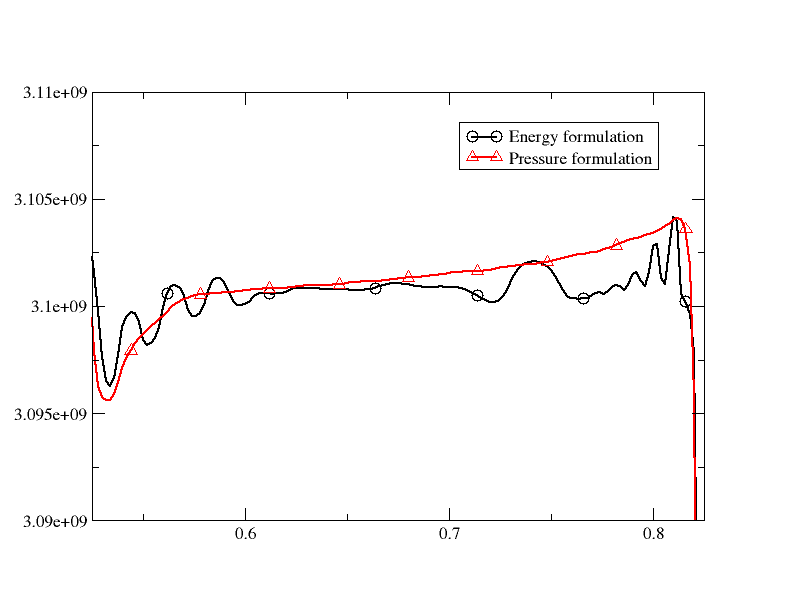}}
\caption{\label{Riemann:cochran} Riemann problem for  the Euler system with Cochran-Chan EOS. Comparison of the solutions obtained from the pressure and energy formulations.
On the right, zooms of the solutions around the contact are provided.}
\end{center}
\end{figure}

The solution displayed in Figure \ref{Riemann:cochran} at a final time $50\cdot 10^{-6}{s}$ show an excellent approximation of the contact discontinuity wave for both nonconservative approximations. It is interesting to notice how the pressure formulation is less oscillatory then the energy one. In general, for both approximations, the shock propagates at the same speed.

\section{Extension to multiphase flows}\label{multiphase}

\subsection{Kapila's five equation model}
Let us consider a different set of equations in the framework of compressible multiphase flows given by the five equation model of Kapila et al. \citep{kapila} and shown in \cite{murrone} to be the formal limit of the Baer and Nunziato model \cite{BN} when the relaxation parameters
simultaneously tend to infinity though being proportional.

\begin{subequations}\label{eq:for:conservative}
\begin{equation}\label{eq:for:conservative:1}
\dpar{\alpha_1}{t}+u\cdot \nabla \alpha_1= K \text{div }u, \qquad K:=\dfrac{\rho_2c_2^2-\rho_1c_1^2}{\frac{\rho_1c_1^2}{\alpha_1}+\frac{\rho_2c_2^2}{\alpha_2}}
\end{equation}
\begin{equation}\label{eq:for:conservative:2}
\dpar{(\alpha_1\rho_1)}{t}+\text{ div }(\alpha_1\rho_1 u)=0
\end{equation}
\begin{equation}\label{eq:for:conservative:3}
\dpar{(\alpha_2\rho_2)}{t}+\text{ div }(\alpha_2\rho_2 u)=0
\end{equation}
\begin{equation}\label{eq:for:conservative:4-2}
\dpar{(\rho u)}{t}+\text{ div }(\rho u^2 + p)=0
\end{equation}
\begin{equation}\label{eq:for:conservative:4}
\dpar{E}{t}+\text{ div }\bigg (( E+p)u\bigg )=0
\end{equation}
\end{subequations}
In this two-phase system, $\alpha_1$ is the volume fraction of phase $\Sigma_1$, while the volume of the second phase $\Sigma_2$ is given by $\alpha_2=1-\alpha_1$. The density of the phase $\Sigma_1$ (respectively $\Sigma_2$) is $\rho_1$ (respectively $\rho_2$). The mixture density is given by $\rho=\alpha_1\rho_1+\alpha_2\rho_2$. We assume a single velocity $u$ and a single pressure $p$. This allows to consider mixture quantities for the energy and momentum conservation laws, while the mass conservation law is described for each phase separately. The internal energy of each phase $\Sigma_j$ is given by $e_j=e_j(p, \rho_j)$ with $j=1,2$, and the mixture internal energy reads $e=\alpha_1e_1+\alpha_2e_2$. The total energy $E$ is the sum of the internal energy and the kinetic energy.
In \eqref{eq:for:conservative:1}, $c_j$ represents the speed of sound of phase $\Sigma_j$ and, in general, this transport equation is nonconservative. System \eqref{eq:for:conservative} is a hyperbolic model and the mixture speed of sound $c$ is defined via
\begin{equation}\label{5eq:son}
\frac{1}{\rho c^2}=\frac{\alpha_1}{\rho_1 c_1^2}+\frac{\alpha_2}{\rho_2 c_2^2}.
\end{equation}
with $c_j$ being the speed of sound of a phase $j$.

The Baer and Nunziato seven equation model \cite{BN}, from which \eqref{eq:for:conservative} has been derived, considers two phases which are described by a set of a mass, momentum and energy conservation laws for each phase and an additional transport equation which links the two sets of equations in terms of the volume fractions. In case of mechanical relaxation, which means that we assume a very large interface between the two phases, we can consider the pressures of each phase to be identical. The same also holds for the velocities.
Following \cite{murrone}, since the pressures of the two phases are equal, their Lagrangian derivatives are equal, too. Therefore, it is possible to write that the entropies $s_1$,$s_2$ are constant and that $p_1(\rho_1,s_1)=p_2(\rho_2,s_2)$, leading to
$$c_2^2\dfrac{d\rho_2}{dt}=c_1^2\dfrac{d\rho_1}{dt}=\frac{\rho_1c_1^2}{\alpha_1} \dfrac{d\alpha_1}{dt}+\rho_1c_1^2\dpar{u}{x}=\frac{\rho_2c_2^2}{\alpha_2} \dfrac{d\alpha_2}{dt}+\rho_2c_2^2\dpar{u}{x}$$

\noindent This allows to reformulate the transport equation as
$$\dfrac{d\alpha_1}{dt}=\frac{\rho_2c_2^2-\rho_1c_1^2}{\frac{\rho_1c_1^2}{\alpha_1}+\frac{\rho_2c_2^2}{\alpha_2}}.$$ and to reduce the original seven equations model of \cite{BN} to the one of \cite{kapila}.

In the following, in order to be able to rewrite the system \eqref{eq:for:conservative} in terms of primitive variables,  we need the differential relations linking the pressure and the internal energy to the densities, $\alpha_j\rho_j$ and the volume fraction $\alpha_1$, since these are independent parameters.
To achieve this, we start from:
$$dp_j=\kappa_j de_j +\chi_j d\rho_j, \qquad \kappa_j=\dpar{p_j}{e_j}\Big{|}_{\rho_j}, \quad  \chi_j=\dpar{p_j}{\rho_j}\Big{|}_{e_j}.$$
Since $e=\alpha_1e_1+\alpha_2e_2$,  $p_1=p_2=p$ and $d(\alpha_j \rho_j)=\alpha_jd\rho_j+\rho_j d\alpha_j$, we get:
\begin{equation*}
\begin{split}
de&=\sum_j\alpha_j d e_j +\sum_j e_jd\alpha_j\\
& =\bigg (\sum_j \frac{\alpha_j}{\kappa_j}\bigg ) dp -\sum_j \dfrac{\chi_j}{\kappa_j}\alpha_jd\rho_j +\sum_j e_jd\alpha_j\\
& =\bigg (\sum_j \frac{\alpha_j}{\kappa_j}\bigg ) dp-\sum_j \dfrac{\chi_j}{\kappa_j}\big ( d(\alpha_j \rho_j)-\rho_j d\alpha_j )+\sum_j e_jd\alpha_j\\
&=\bigg (\sum_j \frac{\alpha_j}{\kappa_j}\bigg ) dp -\sum_j \dfrac{\chi_j}{\kappa_j}d(\alpha_j \rho_j)+ \sum_j \big ( e_j+\rho_j\dfrac{\chi_j}{\kappa_j }\big ) d\alpha_j.
\end{split}
\end{equation*}
We rewrite 
$$e_j+\rho_j \dfrac{\chi_j}{\kappa_j}=\frac{\rho_j }{\kappa_j}\big ( \kappa_j \dfrac{e_j+p}{\rho_j}+\chi_j\big ) -p=\frac{\rho_j{c}_j^2}{\kappa_j} -p,$$
\begin{equation}\label{multiphase:de}
dp=\kappa de+\sum_j \dpar{p}{(\alpha_j\rho_j)} d{(\alpha_j \rho_j)}+\dpar{p}{\alpha_1} d\alpha_1
\end{equation}
with
\begin{equation}
\begin{split}\label{multiphase:derivees}
\frac{1}{\kappa}&=\sum_j \frac{\alpha_j}{\kappa_j}\\
\dpar{p}{(\alpha_j\rho_j)}&=\chi_j \frac{\kappa}{\kappa_j}\\
\dpar{p}{\alpha_1}&=-\kappa\bigg( \frac{\rho_1{c}_1^2}{\kappa_1}-\frac{\rho_2{c}_2^2}{\kappa_2}+e_1-e_2\bigg )
\end{split}
\end{equation}

\subsection{Numerical results}
Similarly as in the previous section, we test the system \eqref{eq:for:conservative} on two different test cases with different equations of state and physical properties. 
\subsubsection{Stiffened Gas EOS}
\label{stiffMF}
To test the capabilities of the nonconservative approximation in the context of compressible multiphase flows, a very severe benchmark problem with high differences in the pressure along a shock tube for epoxy and spinel has been taken from literature \cite{marsh1980,Petitpas2007,Saurel2007,Petitpas2009}.

The considered Riemann problem has a domain of $1m$ length and a discontinuity at $x_d=0.6m$. The parameters for each phase are shown in Table \ref{tab1} and we use the stiffened gas EOS for both phases, which is part of the Mie Gr\"uneisen EOS family and reads 
$p=(\gamma-1)\rho_j e_j -\gamma p_\infty$.\\
{\sout{On the left and on the right of the discontinuity, a} A} mixture of epoxy and spinel is set up, with $u = u_1 = u_2 = 0\; [\frac{m}{s}]$ both on the left and on the right {of the discontinuity}, while on the left the pressure is {set to } $p_L = p_{1L} = p_{2L} = 2 \cdot 10^{11}{[Pa]}$ and on the right {to} $p_R = p_{1R} = p_{2R} = 1\cdot 10^5[Pa]$.
\begin{table}[H]

\begin{center}
\begin{tabular}{llllll}
Phase & Fluid & $\alpha_1$ & $\rho \;[\frac{kg}{m^3}]$ & $\gamma$ & $P_{\infty}\;[Pa]$\\
\hline

1 & Epoxy & $0.5954$ & $1185$ & $2.43$ & $5.3\cdot 10^{9} $\\[0.2em]
2 & Spinel & $0.4046$ & $3622$ & $1.62$ & $141\cdot 10^{9}$\\[0.2em]
\end{tabular}
\end{center}
\caption{\label{tab1} Initial fluid properties and parameters for the test case.}
\end{table}

\begin{figure}[H]
\subfigure[Pressure]{\includegraphics[width=0.45\textwidth]{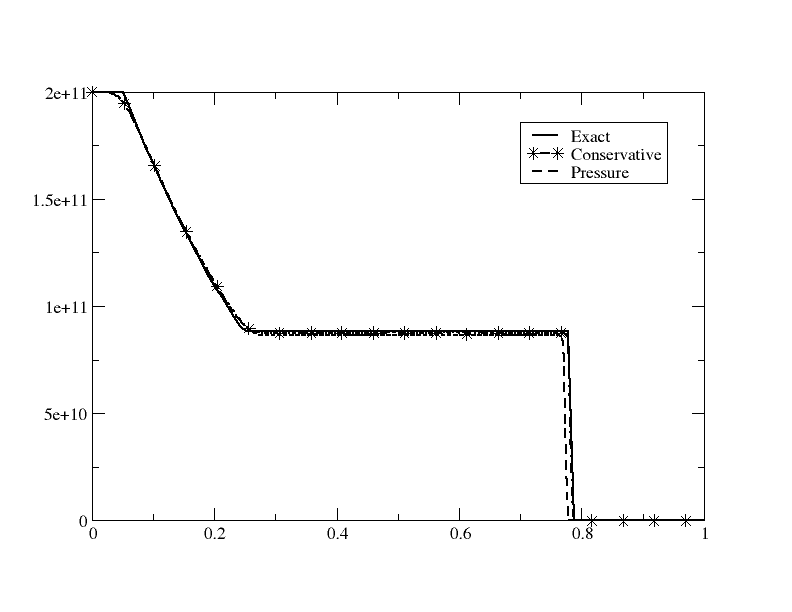}}
\subfigure[Density]{\includegraphics[width=0.45\textwidth]{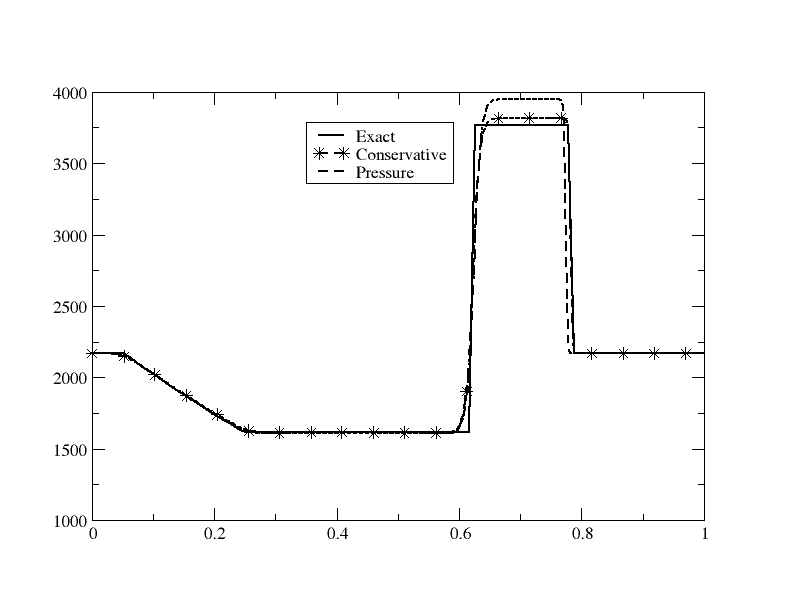}}
\subfigure[Velocity]{\includegraphics[width=0.45\textwidth]{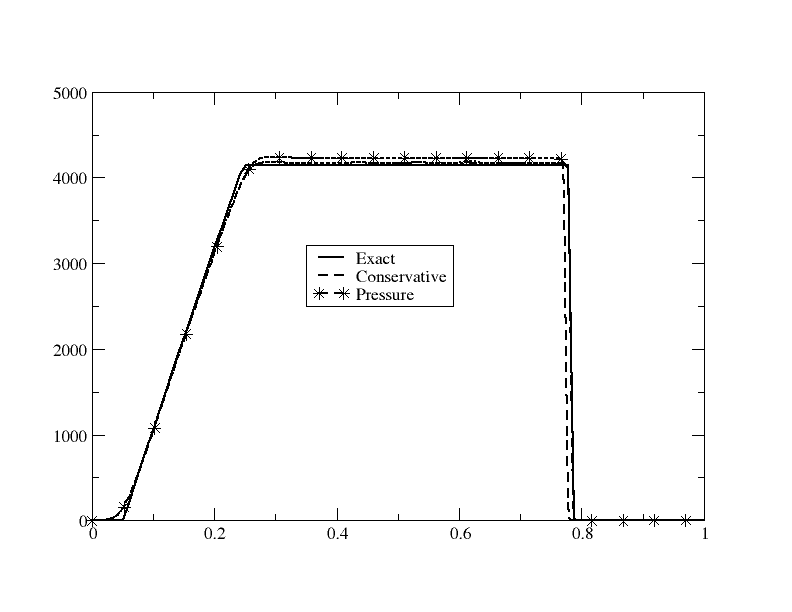}}
\subfigure[$\alpha_1$]{\includegraphics[width=0.45\textwidth]{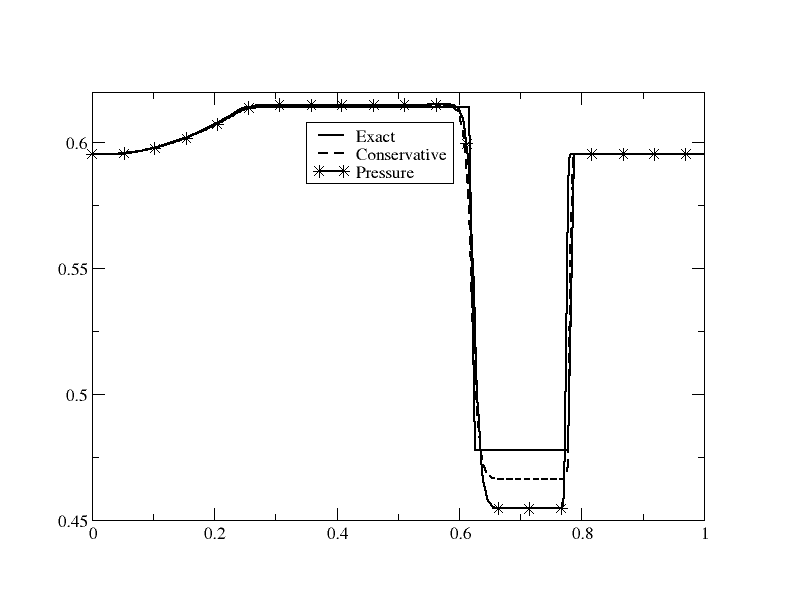}}
\caption{\label{SG} Comparison between the exact solution and numerical solution given by the {classical conservative total energy and nonconservative pressure formulation \sout{schemes}} for the stiffened gas EOS.}
\end{figure}

\begin{figure}[H]
\subfigure[Pressure]{\includegraphics[width=0.45\textwidth]{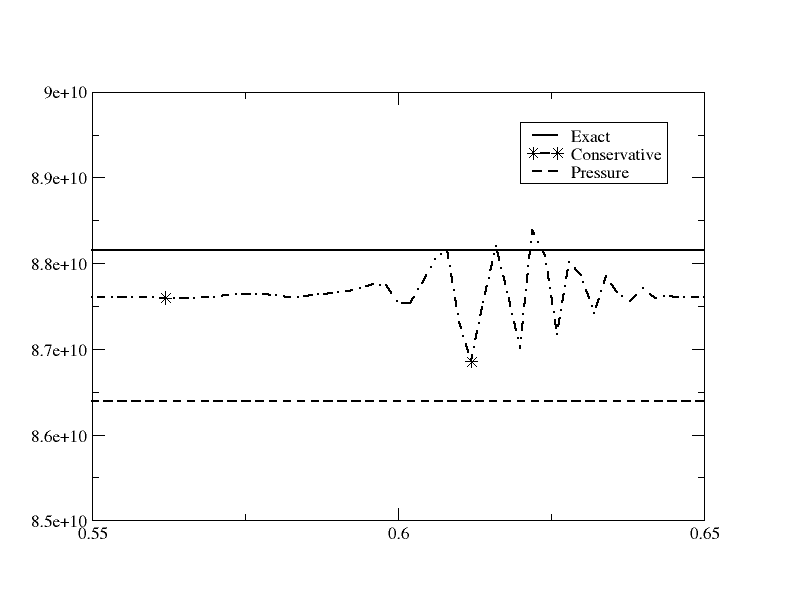}}
\subfigure[Velocity]{\includegraphics[width=0.45\textwidth]{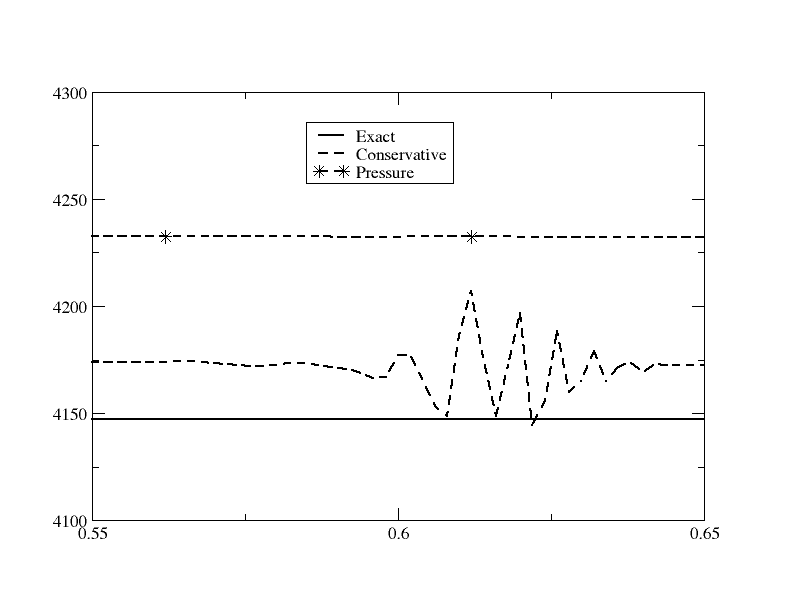}}
\caption{\label{SG:2} Comparison between the exact solution and the numerical solution given by the {classical conservative total energy and nonconservative pressure formulation \sout{schemes}} for the stiffened gas EOS. Zoom of the region of the contact discontinuity.}
\end{figure}

\noindent The final time is $t=29\;\mu s$ and the CFL number is set to $0.5$.
The results are displayed in Figures \ref{SG} and \ref{SG:2}.

The obtained solutions show an excellent contact discontinuity approximation. The difference in the plateaus is caused by the nonconservative form of the systems, which affects the numerical dissipation, while the rarefaction waves are identical. Note however that the complexity of the pressure formulation is lower: the pressure is a primary parameter, and it doesn't need to be computed from the internal energy. In other words, there is no need to get $p$ knowing the mass fraction and the densities via the formula
$$\varepsilon=Y_1\varepsilon_1(\rho_1,p)+ Y_2\varepsilon_2(\rho_2,p).$$
The inversion of the latter relation can be costly for nonlinear EOS. Our goal is certainly not to get the best possible solution in this context, and there exist methods that provide better results on this kind of problems, see e.~g. \cite{saurelNiki,Kumar}. Our goal is only  to show the versatility of our approach.

\subsubsection{Mixed stiffened gas and Cochran-Chan  EOS}
\begin{figure}[H]
\subfigure[Pressure]{\includegraphics[width=0.45\textwidth]{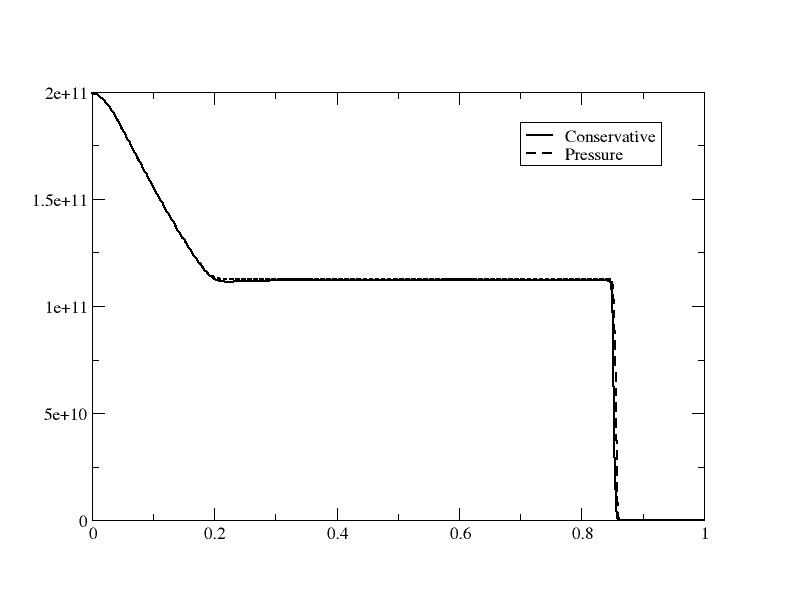}}
\subfigure[Density]{\includegraphics[width=0.45\textwidth]{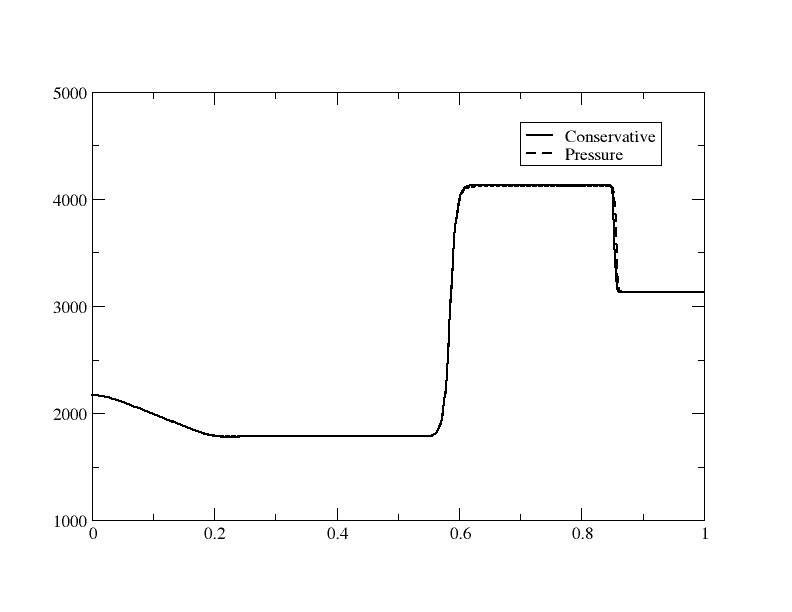}}
\subfigure[Velocity]{\includegraphics[width=0.45\textwidth]{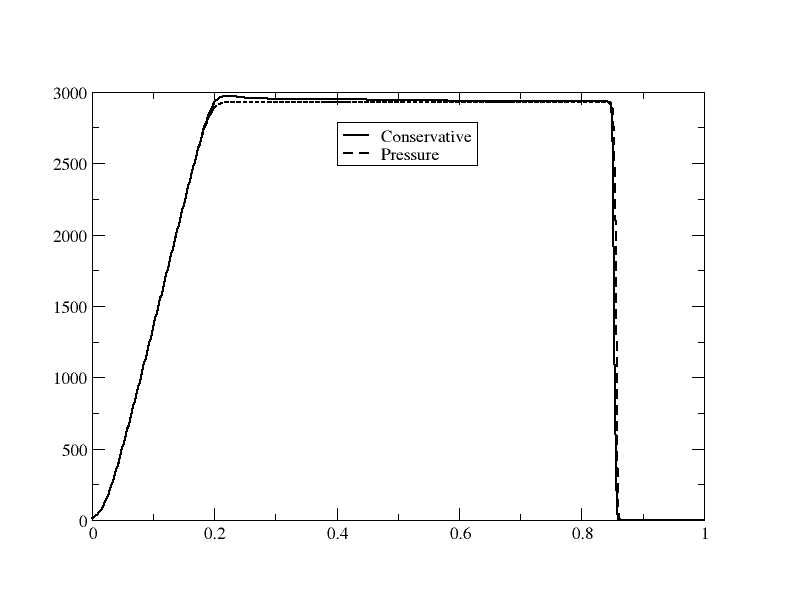}}
\subfigure[$Y_1$]{\includegraphics[width=0.45\textwidth]{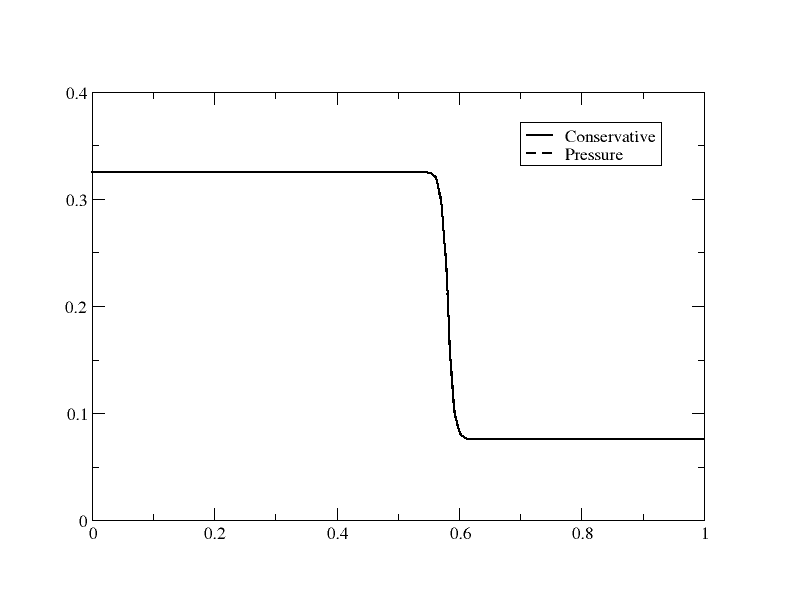}}
\centerline{\subfigure[$\alpha_1$]{\includegraphics[width=0.45\textwidth]{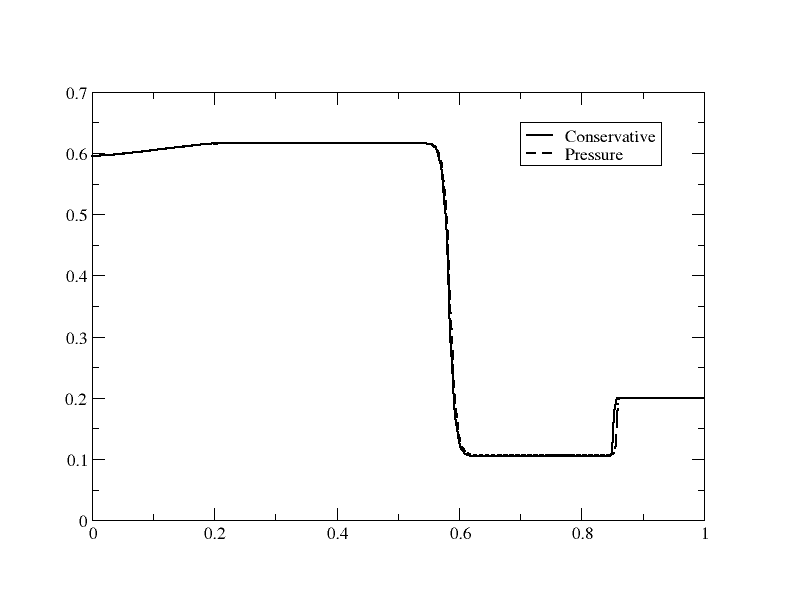}}}
\caption{\label{SG-Mie} Comparison between the exact solution and the numerical solution given by the {classical conservative total energy and nonconservative pressure formulation \sout{schemes}} for a mixture of stiffened gas and Cochran-Chan EOS.}
\end{figure}

\begin{figure}[H]
\subfigure[Pressure]{\includegraphics[width=0.45\textwidth]{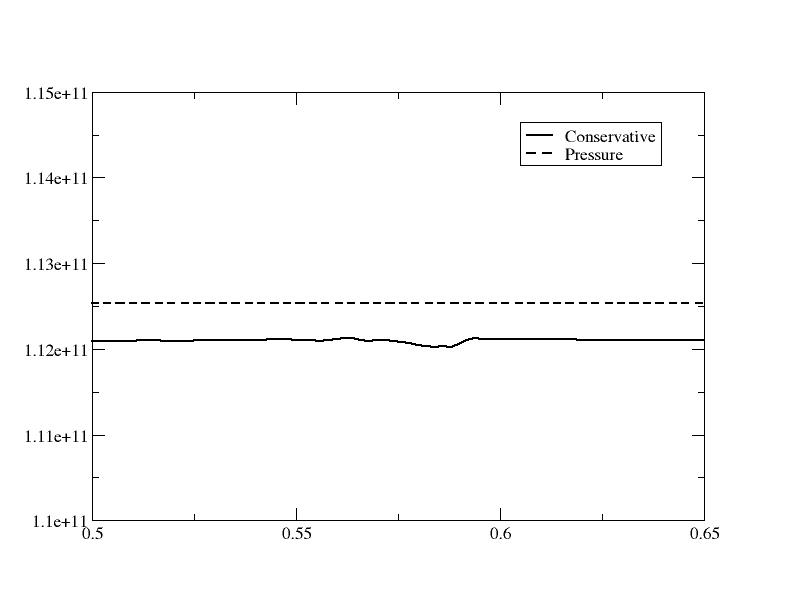}}
\subfigure[Velocity]{\includegraphics[width=0.45\textwidth]{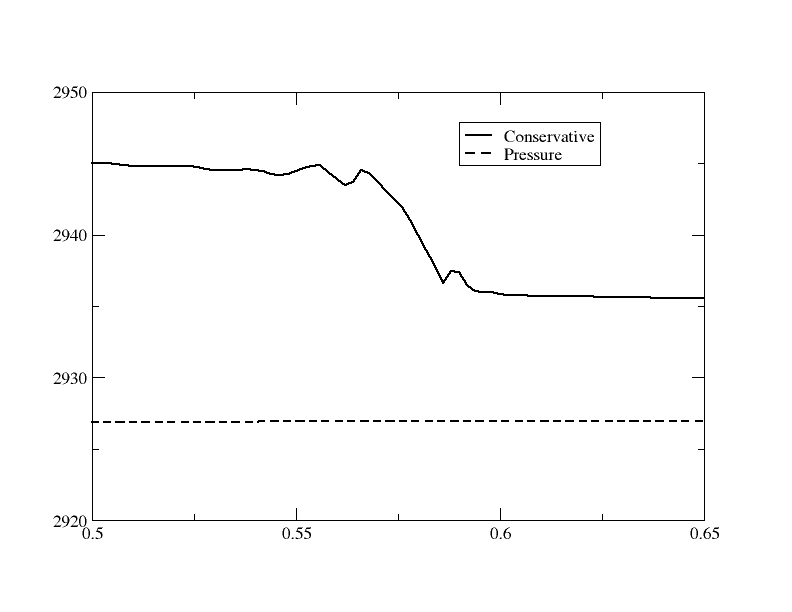}}
\caption{\label{SG-Mie-zoom} Comparison between the exact solution and the numerical solution given by the  {classical conservative total energy and nonconservative pressure formulation \sout{schemes}} for a mixture of stiffened gas and Cochran-Chan EOS. Zoom of the contact discontinuity region.}
\end{figure}

The following test case is a more challenging variation of the previous test case of section \ref{stiffMF}. Let us assume the first phase {to} be described by a stiffened gas EOS with parameters $\gamma_1=2.43$, $p_\infty=5.3\; 10^9$, and the other one by the Cochran-Chan EOS with the parameters  given in Table \ref{tab:cochran}.
We initially set the same conditions as in section \ref{stiffMF} and change only the volume fractions of the involved phases: $\alpha_l=0.5954$ and $\alpha_r=0.2$.
The results are displayed in Figures \ref{SG-Mie} and \ref{SG-Mie-zoom}.

The obtained solutions show again an excellent contact discontinuity approximation without any oscillations. The fact that the plateaus are on different levels is not unexpected, as pointed out before in \ref{stiffMF}.

\section{Conclusions}
\label{Conclusions}
In this paper, we have described a technique that enables the usage of a nonconservative  formulation which nevertheless guarantees the conservation of the involved quantities. This method uses a residual distribution discretization with simple additional conditions which are able to guarantee the convergence to the correct weak solution.  The emphasis is put on nonlinear equations of states where the pressure depends nonlinearly on the density.  The presented approximation is then generalized to a multiphase system. 
The numerical tests have been done in one dimension in order to compare the results to the exact solutions. 
Extension to multidimensional problems  can be done following the lines of \cite{ricchiuto} for the second order.
Extension to higher order {of accuracy} could easily be obtained using the ideas of \cite{abgrall-Shu,abgrallPaolaTokareva}, where the Runge-Kutta type timestepping scheme that we use here can be reinterpreted as a particular version of a deferred correction method. In this case, the relation \eqref{cons:e total} stays the same provided the residuals are defined according to \cite{abgrall-Shu} and the residual on the internal energy behaves like $O(h^{k+d})$, where $k$ is the expected order and $d$ the spatial dimension. We emphasize that the algebra remains identical. This will be demonstrated in a forthcoming paper.

\section*{Acknowledgements.}
P.~B. and S.~T. have been funded by SNSF project \# 200021\_153604. R.~A. has been funded in part by the same grant.

\appendix
\section{A remark on the convergence to a weak solution}\label{appendix A}
There are two ways of showing this. We provide the main idea for the system \eqref{eq:nc}.
Assuming that we have a scheme for this system that satisfies \eqref{cons:e total}, and to simplify the derivations, we shall use the first order time scheme.
Then one can define a residual for the total energy by simply setting
$${\Phi_{\sigma,E}^{K}}:={\Phi_{\sigma,e}^{K}}+\frac{\overline{\mathbf{u}_\sigma}+\underline{\mathbf{u}_\sigma} }{2}\cdot {\Phi_{\sigma,\textbf{m}}^{K}}-\frac{1}{2}\overline{\bu}_\sigma\cdot \underline{\bu}_\sigma\;{\cdot \, \Phi_{\sigma,\rho}^{K}}.
$$
By construction we have
$$\sum_{\sigma\in K}{\Phi_{\sigma,E}^{K}}=\Phi_E^K.$$
This shows that the sequence of solutions will converge to a weak solution. Using the results of \cite{icm}, one can compute numerical fluxes for the density, momentum and total energy, so that local conservation is guaranteed. It is also simple to extend the proof of the Lax-Wendroff type theorem of \cite{Abgrall:Roe} using the same conditions.

\section*{References}
\bibliography{biblio}
\bibliographystyle{elsarticle-num}
\biboptions{sort&compress}

\end{document}